\documentclass{amsart}

\usepackage[mathscr]{eucal}
\usepackage{amsmath}
\usepackage{amssymb}
\usepackage{amscd}
\usepackage{amsfonts}
\usepackage{enumerate}

\usepackage[all]{xy} 
\CompileMatrices

\theoremstyle{plain}
\newtheorem{thm}{Donotwrite}[section]

\newtheorem{define}[thm]{Definition}
\newtheorem{theorem}[thm]{Theorem}
\newtheorem{proposition}[thm]{Proposition}
\newtheorem{lemma}[thm]{Lemma}

\theoremstyle{definition}
\newtheorem{example}[thm]{Example}

\newtheorem{remark}[thm]{Remark}

\numberwithin{equation}{section}

\newfont{\germ}{eufm10}

\newcommand\cd{\cdots}

\newcommand\et[1]{\tilde{e}_{#1}}
\newcommand\ft[1]{\tilde{f}_{#1}}
\newcommand\geh{\goth{g}}
\newcommand\gehl[1]{\geh_{#1}}
\newcommand\goth[1]{\mbox{\germ #1}}
\newcommand\Kz{K_{\Z}}
\newcommand\La{\Lambda}
\newcommand\la{\lambda}
\newcommand\lw[1]{\lower.4mm\hbox{${#1}$}}
\newcommand\ol{\overline}
\newcommand\olLa{\ol{\La}}
\newcommand\ot{\otimes}
\newcommand\Q{{\mathbb Q}}
\newcommand\veps{\varepsilon}
\newcommand\vphi{\varphi}
\newcommand\wt{\mbox{\sl wt}\,}
\newcommand\Z{{\mathbb Z}}
\newcommand\Zn{\Z_{\ge0}}
\newcommand{\bl}{\bigl}
\newcommand{\br}{\bigr}
\newcommand{\seteq}{\mathbin{:=}}
\newcommand{\set}[2]{\left\{{#1}\ \vert\ {#2}\right\}}
\newcommand{\be}{\begin{enumerate}}
\newcommand{\ee}{\end{enumerate}}
\newcommand{\cl}{\colon}
\newcommand{\hs}{\hspace*}
\newcommand{\isoto}[1][]{\mathop{\xrightarrow[#1]%
{\rule{0pt}{.9ex}%
{\raisebox{-.4ex}[0ex][-.65ex]{$\mspace{1mu}\sim\mspace{2mu}$}}}}}
\newcommand{\g}{\mathfrak{g}}
\newcommand{\eq}{\begin{eqnarray}}
\newcommand{\eneq}{\end{eqnarray}}
\newcommand{\Lemma}{\begin{lemma}}
\newcommand{\enlemma}{\end{lemma}}
\newcommand{\ba}{\begin{array}}
\newcommand{\ea}{\end{array}}

\title{Perfect Crystals for $U_q(D_4^{(3)})$} 

\author{M. Kashiwara}
\address{Masaki Kashiwara: 
Research Institute for Mathematical Sciences, Kyoto University, Kyoto 606-8502, Japan}

\author{K.C. Misra}
\address{Kailash C. Misra: 
Department of Mathematics, North Carolina State University, Raleigh, North Carolina 27695-8205, USA} 

\author{M. Okado}
\address{Masato Okado: 
Department of Mathematical Science, Graduate School of Engineering Science, 
Osaka University, Toyonaka, Osaka 560-8531, Japan}

\author{D. Yamada}
\address{Daisuke Yamada: 
Graduate School of Mathematical Sciences, the University of Tokyo, 
Megro, Tokyo 5465-8274, Japan} 

\date{}

\begin{document}

\maketitle

\section{Introduction}

In \cite{KMN1} the notion of {\em perfect crystal} was introduced. In the subsequent paper
\cite{KMN2} a perfect crystal of an arbitrary level was given for the quantum affine algebra
corresponding to every non-exceptional affine algebra. Later studies revealed 
that there are more perfect crystals than listed in \cite{KMN2}. See \cite{BFKL,Ko,JMO,SS,Y,Ka3,NS}
for example. Actually, there is a conjecture originating from {\em fermionic formulas} 
\cite{HKOTY,HKOTT}, saying that a certain finite-dimensional module, called Kirillov-Reshetikhin 
module, KR module for short, of a quantum affine algebra has a crystal base. KR modules are 
parametrized by two integers $r,s$. $r$ corresponds to a vertex of the Dynkin diagram of the affine
algebra except a distinguished vertex 0 as in \cite{Kac} and $s$ is a positive integer. The conjecture
also states that if $s$ is a multiple of $t_r\seteq
\max(1,2/(\alpha_r,\alpha_r))$, then the conjectural 
crystal base of the corresponding KR module is perfect of level $s/t_r$. Here $(\cdot,\cdot)$ stands
for the standard bilinear form on the weight lattice as in \cite{Kac}. The latter conjecture was 
derived by taking a suitable limit of the corresponding fermionic formula.
Up to now there is no counterexample to this conjecture. We remark that apart from the existence 
of crystal base, KR modules have marvelous properties, such as algebraic relations among characters
\cite{KR,KNS,N,H}, relation with Demazure modules \cite{KMOTU,K4,FL,CM,FSS}, that with fusion products 
\cite{CL,AK} and so on.

In this paper, we add another evidence for the conjecture to be valid. To explain what we have done
more precisely let us recall basic notations and the definition of perfect crystal. 
Let $\mathfrak{g}$ be an affine algebra, $U_q(\mathfrak{g})$ the associated quantum affine algebra 
and $U'_q(\mathfrak{g})$ its subalgebra without the degree operator.
Let $P$ be the weight lattice, i.e., $P=\sum_i\Z\La_i\oplus\Z\delta$ where $\La_i$ is a fundamental
weight and $\delta$ is the generator of the null roots, and set 
$P_{cl}=P/\Z\delta$,
$P_{cl}^+=\set{\lambda \in P_{cl}}%
{\langle h_i, \lambda \rangle \geq 0 \ \text{for any $i$}}
=\sum_i\Z_{\ge0}\Lambda_i$
and $(P_{cl}^+)_l=\{\lambda \in P_{cl}^+\vert
\langle c, \lambda \rangle =l \}$ for $l\in\mathbb{Z}_{\geq 0}$. Here $c$ is the canonical central 
element. Let $\mathrm{Mod}^f(\mathfrak{g},P_{cl})$ be the category of finite-dimensional 
$U'_q(\mathfrak{g})$-modules.
The modules in this category 
have the weight decomposition with respect to $P_{cl}$.
Let $B$ be a $P_{cl}$-weighted crystal. For $b\in B$, we set 
$\varepsilon(b)=\sum_i\varepsilon_i(b)\La_i$ and $\varphi(b)=\sum_i\varphi_i(b)\La_i$, where
$\veps_i(b)=\max\{k\mid\et{i}^kb\neq0\}$ and $\vphi_i(b)=\max\{k\mid\ft{i}^kb\neq0\}$.

\begin{define} \label{def:perfect}
$([\mathrm{KMN1}])$
For $l\in \mathbb{Z}_{>0}$ we say $B$ is a perfect crystal of 
level $l$ if $B$ satisfies the following conditions. 
\begin{itemize}
\item[$(P1)$] $B\otimes B$ is connected. 
\item[$(P2)$] There exists $\lambda_0\in P_{cl}$ such that 
$\mathrm{wt}(B)\subset \lambda_0+\sum_{i\neq0}\mathbb{Z}_{\leq 0}\alpha_i$ 
and that $\sharp (B_{\lambda_0})=1$. 
\item[$(P3)$] There is a $U'_q(\mathfrak{g})$-module in 
$\mathrm{Mod}^f(\mathfrak{g},P_{cl})$ with a crystal pseudobase $(L,B')$ 
such that $B$ is isomorphic to $B'/\{\pm 1\}$. 
\item[$(P4)$] For any $b\in B$, we have 
$\langle c, \varepsilon(b)\rangle\geq l$.
\item[$(P5)$] The maps $\varepsilon$ and $\varphi$ from 
$B_{\min}\seteq\set{b\in B}{\langle c, \varepsilon(b)\rangle =l}$ to 
$(P_{cl}^+)_l$ are bijective. 
\end{itemize}
\end{define}

We consider the quantum affine algebra $U'_q(\mathfrak{g})$ corresponding to the exceptional
affine algebra $\mathfrak{g}=D_4^{(3)}$. The KR modules we treat in this paper correspond to the 
vertex $1$ in the Dynkin diagram in section 2.1. According to the above conjecture the KR module 
for the pair $(r,s)=(1,l)$ is to be perfect of level $l$. Let us look through the main text in order.
In section 2 we construct a fundamental representation $V^1$ of $U_q(D_4^{(3)})$ and calculate the
quantum $R$-matrix that will be used in the next section. In section 3 we construct a family of 
$U'_q(D_4^{(3)})$-modules $\{V^l\}_{l\ge1}$ by so-called fusion construction. This $V^l$ is nothing
but the KR module for $(1,l)$. Using a technique developed in \cite{KMN2} one sees that $V^l$
has a crystal pseudobase (Theorem \ref{th:V^l}), which confirms (P2) and (P3) of Definition 
\ref{def:perfect}. In section 4 we introduce a $U_q(G_2)$-crystal $B_l$ that has the same 
decomposition as $V^l$ and define the $0$-action by hand. It was obtained by the investigation of
the result of \cite{Y} and computer experiments.
We check in section 5 that with respect to the $(0,1)$-actions $B_l$ turns into a disjoint 
union of $U_q(A_2)$-crystals (Theorem \ref{thm:$U_q(A_2)$CrystalStructure}). Hence by Theorem 
\ref{th:uniqueness} that states the uniqueness of such crystal, one concludes that $B_l$ is the 
crystal base of $V^l$. We stress here that this theorem can be applied to affine algebras of other 
types. (P1), (P4) and (P5) are checked in the last section. Thus we have obtained 

\begin{theorem}
For $l\in \mathbb{Z}_{>0}$
$B_l$ is a perfect crystal of level $l$. 
\end{theorem}
The value of this theorem lies in the following fact. Let $B(\la)$ denote the crystal of the 
irreducible highest weight $U_q(\mathfrak{g})$-module of a dominant integral weight $\la$. If $B_l$
is a perfect crystal of level $l$, then for any dominant integral weight $\la$ of level $l$ there
exists a unique dominant integral weight $\mu$ and an isomorphism of crystals
\[
B(\la)\simeq B(\mu)\ot B_l.
\]
By iterating it, one can obtain the so-called Kyoto path model of $B(\la)$. By analogy one may 
consider a path model for the crystal $B(\infty)$ of the negative part $U_q^-(\mathfrak{g})$ of 
$U_q(\mathfrak{g})$. For this purpose one needs a coherent family of perfect crystals (Definition
\ref{def:coherent}). At the final stage we show that our family of perfect crystals $\{B_l\}_{l\ge1}$
are coherent. It implies that there is a limit $B_{\infty}$ of $\{B_l\}_{l\ge1}$ and an isomorphism
of crystals
\[
B(\infty)\simeq B(\infty)\ot B_{\infty}.
\]

\section{$U'_q(D_4^{(3)})$ and its fundamental representation}

\subsection{Quantum affine algebra $U'_q(D_4^{(3)})$}

We collect necessary data for the affine Lie algebra $D_4^{(3)}$.
Let 
$\{\alpha_0, \alpha_1, \alpha_2\}$, $\{h_0, h_1, h_2\}$ and 
$\{\La_0, \La_1, \La_2\}$ be the set of 
simple roots, simple coroots and fundamental weights, respectively.
The generalized Cartan matrix $(\langle h_i,\alpha_j\rangle)_{i,j=0,1,2}$
is given by
\[
\left(
\begin{array}{rrr}
2  & -1 & 0  \\
-1 & 2  & -3 \\
0  & -1 & 2
\end{array}
\right),
\]
and its Dynkin diagram is depicted as follows.
\[\xymatrix@R=.0ex{
0&1&2\\
\circ\ar@{-}[r]&\circ&\circ\ar@3{->}[l]
}\]
The standard null root $\delta$ and the canonical central element $c$ are 
given by
\[
\delta=\alpha_0+2\alpha_1+\alpha_2\quad\text{and}\quad c=h_0+2h_1+3h_2.
\]
The affine Lie algebra $D_4^{(3)}$ contains 
the finite-dimensional simple Lie algebra $G_2$. Its fundamental
weights are given by $\olLa_1=\La_1-2\La_0$, $\olLa_2=\La_2-3\La_0$.

The quantum affine algebra $U'_q(D_4^{(3)})$ is 
the associative algebra over $\Q(q)$ generated by 
$\{e_i, f_i, t_i^{\pm1} \mid i=0,1,2\}$ satisfying the relations: 
\begin{align*}
&t_i t_j=t_j t_i, \quad 
t_i t_i^{-1}=t_i^{-1}t_i=1, 
\\
&t_i e_j t_i^{-1}=q_i^{\langle h_i,\alpha_j \rangle} e_j, \quad 
t_i f_j t_i^{-1}=q_i^{-\langle h_i, \alpha_j \rangle} f_j, \quad 
[e_i, f_j]=\delta_{ij}\frac{t_i-t_i^{-1}}{q_i-q_i^{-1}}, 
\\
&\sum_{n=0}^l (-1)^n e_i^{(n)} e_j e_i^{(l-n)}
=\sum_{n=0}^l (-1)^n f_i^{(n)} f_j f_i^{(l-n)}=0
\;
\text{where $i\neq j$, $l=1-\langle h_i, \alpha_j \rangle$.}
\end{align*}
Here we use the following notations: $[m]_i=(q_i^m-q_i^{-m})/(q_i-q_i^{-1})$,
$[n]_i!=\prod_{m=1}^n[m]_i$, $e_i^{(n)}=e_i^n/[n]_i!,f_i^{(n)}=f_i^n/[n]_i!$ with
\[
q_0=q_1=q,\quad q_2=q^3.
\]
The algebra $U_q(D_4^{(3)})$ is defined by introducing another generator $q^d$, so 
$U'_q(D_4^{(3)})$ is its subalgebra. 
The algebra $U'_q(D_4^{(3)})$ has subalgebras 
$U_q(G_2)$ and $U_q(A_2)$ generated 
by $\{e_i,f_i,t_i^{\pm1}\mid i=1,2\}$
and $\{e_i,f_i,t_i^{\pm1}\mid i=0,1\}$, respectively.

\subsection{Fundamental representation} \label{subsec:fund rep}

Let $V^1=V^{G_2}(\olLa_1)\oplus V^{G_2}(0)$ denote the direct sum of 
the irreducible highest weight $U_q(G_2)$-modules 
$V^{G_2}(\olLa_1)$ and $ V^{G_2}(0)$ with highest weights $\olLa_1$ 
and $0$. Let $\{ v_1, v_2, v_3, v_0, v_{{\bar 3}}, v_{{\bar 2}}, v_{{\bar 1}}, v_{\phi}\}$ 
be a basis of $V^1$, where $v_{\phi}$ belongs to $V^{G_2}(0)$. $V^1$ is endowed with the 
$U'_q(D_4^{(3)})$-module structure. First, if $v$ is a weight vector of weight $\la$,
we have $t_iv=q_i^{\langle h_i,\la\rangle}v$. 
The weight of each basis vector is given
respectively by $\La_1-2\La_0,\La_2-\La_1-\La_0,2\La_1-\La_2-\La_0,0,\La_0+\La_2-2\La_1,
\La_0+\La_1-\La_2,2\La_0-\La_1,0$. The actions of $e_i,f_i$ are given as follows.
\begin{align*}
&e_0 v_1=v_{\phi}+\frac{1}{[2]}v_0, \ 
 e_0 v_2=v_{{\bar 3}}, \
 e_0 v_3=v_{{\bar 2}}, \
 e_0 v_0=v_{{\bar 1}}, \
 e_0 v_{\phi}=\frac{[3]}{[2]}v_{{\bar 1}}, \\
&f_0 v_{{\bar 1}}=v_{\phi}+\frac{1}{[2]}v_0, \
 f_0 v_{{\bar 2}}=v_3, \
 f_0 v_{{\bar 3}}=v_2, \
 f_0 v_0=v_1, \
 f_0 v_{\phi}=\frac{[3]}{[2]}v_1, \\
&e_1 v_2=v_1, \
 e_1 v_0=[2]v_3, \
 e_1 v_{{\bar 3}}=v_0,  \
 e_1 v_{{\bar 1}}=v_{{\bar 2}}, \\
&f_1 v_{{\bar 2}}=v_{{\bar 1}}, \
 f_1 v_0=[2]v_{{\bar 3}}, \
 f_1 v_3=v_0, \
 f_1 v_1=v_2, \\
&e_2 v_3=v_2, \
 e_2 v_{{\bar 2}}=v_{{\bar 3}}, \\
&f_2 v_{{\bar 3}}=v_{{\bar 2}}, \
 f_2 v_2=v_3.
\end{align*}
If the action of some basis vector is not written, then we should understand
that it is $0$. Hereafter $[m]$ always means $[m]_0=[m]_1$.

Set $A_{\Z}=\{f(q)/g(q)\mid f(q),g(q)\in\Z[q],g(0)=1\}$ and $\Kz=A_{\Z}[q^{-1}]$.
We define $U'_q(\geh)_{\Kz}$ as the $\Kz$-subalgebra of $U'_q(\geh)$ generated by $e_i,f_i,t_i^{\pm1}$.
It is easy to see that $V^1$ admits a $U'_q(\geh)_{\Kz}$-submodule $V^1_{\Kz}$.
In the subsequent section we need a polarization $(\;,\;)$ on $V^1$ such that 
$(V^1_{\Kz},V^1_{\Kz})\subset\Kz$. See section 2 of
\cite{KMN2} for the polarization. It is constructed as follows. It is known (\cite{K1})
that for any dominant integral weight $\la$ 
the irreducible highest weight $U_q(G_2)$-module 
$V^{G_2}(\la)$ of highest weight $\la$ has a polarization. Let $(\;,\;)_1$ be such 
polarization on $V^{G_2}(\olLa_1)$ normalized as $(v_1,v_1)_1=1$. 
We define a symmetric
bilinear form $(\;,\;)$ on $V^1$ by requiring
\begin{align*}
(u,v)&=(u,v)_1\text{ for }u,v\in V^{G_2}(\olLa_1),\\
(u,v_{\phi})&=0\text{ for }u\in V^{G_2}(\olLa_1),\\
(v_{\phi},v_{\phi})&=q\frac{[3]}{[2]}.
\end{align*}
Then $(\;,\;)$ satisfies 
\[
(t_iu,v)=(u,t_iv),\quad (e_iu,v)=(u,q_i^{-1}t_i^{-1}f_iv),\quad 
(f_iu,v)=(u,q_i^{-1}t_ie_iv)
\]
for all $u,v\in V^1$ and it becomes a polarization.
$(V^1_{\Kz},V^1_{\Kz})\subset\Kz$ can also be checked.

As a $U_q(G_2)$-module, the tensor product $V^1\otimes V^1$ decomposes into 
\begin{equation}\label{eq:decomp}
V^1\otimes V^1 \simeq 
V^{G_2}(2\olLa_1)\oplus V^{G_2}(\olLa_2)\oplus V^{G_2}(\olLa_1)^{\oplus 3}
\oplus V^{G_2}(0)^{\oplus 2}. 
\end{equation}
A highest weight vector in each irreducible component is listed below.
\begin{align*}
u_{2\La_1}=&v_1\otimes v_1, 
\\
u_{\La_2}=&v_1\otimes v_2-qv_2\otimes v_1, 
\\
u_{\La_1}^{(1)}=&v_1\otimes v_{\phi}, 
\\
u_{\La_1}^{(2)}=&v_{\phi}\otimes v_1, 
\\
u_{\La_1}^{(3)}=&v_1\otimes v_0-q^6v_0\otimes v_1-q^2[2]v_2\otimes v_3+
q^5[2]v_3\otimes v_2, 
\\
u_0^{(1)}=&v_{\phi}\otimes v_{\phi}, 
\\
u_0^{(2)}=&v_1\otimes v_{{\bar 1}}+q^{10}v_{{\bar 1}}\otimes v_1
-qv_2\otimes v_{{\bar 2}}-q^9v_{{\bar 2}}\otimes v_2
+q^4v_3\otimes v_{{\bar 3}}
\\&
+q^6 v_{{\bar 3}}\otimes v_3-\frac{q^4}{[2]}v_0\otimes v_0.  
\end{align*}
Here lower indices signify the highest weights. For the action on the 
tensor product, we use the lower coproduct. Namely, 
\[
\Delta(e_i)=e_i\ot t_i^{-1}+1\ot e_i,\quad \Delta(f_i)=f_i\ot1+t_i\ot f_i.
\]

\subsection{Calculation of the $R$-matrix}

Let $V_x^1=\Q[x,x^{-1}]\ot V^1$ be the $U'_q(D_4^{(3)})$-module 
with the actions 
of $e_i,f_i,t_i$ replaced with $x^{\delta_{i0}}e_i$, $x^{-\delta_{i0}}f_i$, $t_i$,
respectively. The $R$-matrix $R(x,y)$ for $V^1\ot V^1$ is an operator 
\[
R(x, y) \; : \; V^1_x \otimes V^1_y \longrightarrow V^1_y \otimes V^1_x
\]
commuting with the actions of $U_q'(D_4^{(3)})$. It is unique up to a scalar 
multiple and satisfies the following properties. 
\begin{itemize}
\item[(1)] $R(x,y)\in\Q(q)[x/y,y/x]\ot\text{End}(V^1\ot V^1)$.
\item[(2)] The Yang-Baxter equation holds:
\begin{align*}
(R(y, z)\otimes 1)(1\otimes &R(x, z))(R(x, y)\otimes 1)\\
&=(1\otimes R(x, y))(R(x, z)\otimes 1)(1\otimes R(y,z)).
\end{align*}
\item[(3)] $R(x,y)R(y,x)\in\Q(q)[x/y,y/x]$.
\end{itemize}

By $U_q'(G_2)$-linearity and (\ref{eq:decomp}), we have 
\begin{align*}
R(u_{2\La_1})&=a^{2\La_1}u_{2\La_1},& R(u_{\La_2})&=a^{\La_2}u_{\La_2},\\
R(u_{\La_1}^{(i)})&=\sum_{j=1}^3a^{\La_1}_{ij}u_{\La_1}^{(j)}\;(i=1,2,3),&
R(u_0^{(i)})&=\sum_{j=1}^2a^0_{ij}u_0^{(j)}\;(i=1,2).
\end{align*}
To calculate these coefficients we prepare

\begin{lemma}
We have the following relations.
\begin{itemize}
\item[1.] $f_0f_1f_2u_{\La_2}=(qxy)^{-1}(x-q^2y)u_{2\La_1}$
\item[2.] $f_0u_{\La_1}^{(1)}=(q^2y)^{-1}[3]/[2]u_{2\La_1}$
\item[3.] $f_0u_{\La_1}^{(2)}=x^{-1}[3]/[2]u_{2\La_1}$
\item[4.] $f_0u_{\La_1}^{(3)}=(q^2xy)^{-1}(x-q^8y)u_{2\La_1}$
\item[5.] $f_0^2f_1f_2f_1u_{\La_1}^{(1)}=(qxy)^{-1}[3]u_{2\La_1}$
\item[6.] $f_0^2f_1f_2f_1u_{\La_1}^{(2)}=(qxy)^{-1}[3]u_{2\La_1}$
\item[7.] $f_0^2f_1f_2f_1u_{\La_1}^{(3)}=(qxy)^{-2}[2]([2](x^2-q^8y^2)-q^3(1-q^2)xy)u_{2\La_1}$
\item[8.] $f_0^3(f_1f_2f_1)^2u_{\La_1}^{(1)}=(x^2y)^{-1}[2][3]^2u_{2\La_1}$
\item[9.] $f_0^3(f_1f_2f_1)^2u_{\La_1}^{(2)}=(q^2xy^2)^{-1}[2][3]^2u_{2\La_1}$
\item[10.] $f_0^3(f_1f_2f_1)^2u_{\La_1}^{(3)}=(qxy)^{-2}[2]^2[3](x-q^8y)u_{2\La_1}$
\item[11.] $f_0^2u_0^{(1)}=(qxy)^{-1}[3]^2/[2]u_{2\La_1}$
\item[12.] $f_0^2u_0^{(2)}=(q^2xy)^{-2}([2](x^2+q^{14}y^2)-q^7xy)u_{2\La_1}$
\item[13.] $f_0^3(f_1f_2f_1)^2f_0u_0^{(1)}=(q^2x^3y^3)^{-1}[3]^3(x^2+q^2y^2)u_{2\La_1}$
\item[14.] $f_0^3(f_1f_2f_1)^2f_0u_0^{(2)}
	=(q^4x^3y^3)^{-1}[2][3]((x-q^6y)(x-q^8y)+q^2[3](1+q^{10})xy)u_{2\La_1}$
\end{itemize}
\end{lemma}
{}From this one can calculate the coefficients $a^{2\La_1},a^{\La_2},a^{\La_1}_{ij}
(i,j=1,2,3),a^0_{i,j}(i,j=1,2)$. Let $P_{2\La_1},P_{\La_2},P_{\La_1}^{(i)}(i=1,2,3),
P_0^{(i)}(i=1,2)$ be the projections from $V^1\ot V^1$ onto $U_q(G_2)$-submodule 
$V^{G_2}(2\olLa_1),V^{G_2}(\olLa_2),U_q(G_2)u_{\La_1}^{(i)},U_q(G_2)u_0^{(i)}$, respectively.
Let $\iota_{\La_1}^{(i,j)}(i,j=1,2,3)$ 
(resp. $\iota_0^{(i,j)}(i,j=1,2)$) be the 
$U_q(G_2)$-isomorphism sending $u_{\La_1}^{(j)}$ to $u_{\La_1}^{(i)}$ 
(resp. $u_0^{(j)}$ to $u_0^{(i)}$).
Then we have the spectral decomposition of the $R$-matrix.

\begin{proposition} \label{prop:R mat}
Let $z=x/y$. Up to a multiple of an element of $\Q(q)(z)$, the $R$-matrix is of the 
following form
\begin{align*}
R(x,y)=&(1-q^2z)(1-q^6z)(1+q^4z+q^8z^2)P_{2\La_1}\\
&+(z-q^2)(1-q^6z)(1+q^4z+q^8z^2)P_{\La_2}\\
&+\sum_{i,j=1}^3a^{\La_1}_{ij}\iota^{(j,i)}_{\La_1}P_{\La_1}^{(i)}
+\sum_{i,j=1}^2a^0_{ij}\iota^{(j,i)}_0P_0^{(i)},
\end{align*}
where $a^{\La_1}_{ij},a^0_{ij}$ are given by
\begin{align*}
a^{\La_1}_{11}&=a^{\La_1}_{22}=(1-q^6)z(1-q^{12}z^2)/(1+q^2),\\
a^{\La_1}_{12}&=a^{\La_1}_{21}=q^2(1-z)(1-q^6z)
\bigl((1+q^2)(1+q^6z^2)+(q^2+q^6)z\bigr)/(1+q^2),\\
a^{\La_1}_{13}&=a^{\La_1}_{23}=q(1-q^6)z(1-z)(1-q^6z)/(1+q^2)^2,\\
a^{\La_1}_{31}&=a^{\La_1}_{32}=(1+q^2)^2(1+q^8)a^{\La_1}_{13},\\
a^{\La_1}_{33}&=(1-q^6z)
\bl(q^2(1+q^2)(z^3-q^6)+(1-q^2)(1-q^6)z(z-q^4)\br)/(1+q^2),\\
a^0_{11}&=((1+q^2)(q^2+q^{14}z^4)-(1+q^8)z(q^2+q^8z^2)\\
	&\hspace{20ex}+(1-q^4)(1-q^6)(1+q^8)z^2)/(1+q^2),\\
a^0_{12}&=q(1-q^6)^2z(1-z^2)/\bl((1+q^2)(1-q^4)\br),\\
a^0_{21}&=q(1-q^{14})(1-q^4+q^8)z(1-z^2),\\
a^0_{22}&=z^4(a^0_{11}|_{z\rightarrow1/z}).
\end{align*}
Moreover,
\begin{align*}
\det(a^{\La_1}_{ij})&=\frac{(z-q^2)^2(z^2+q^4z+q^8)}{(1-q^2z)^2(1+q^4z+q^8z^2)}(a^{2\La_1})^3,\\
\det(a^0_{ij})&=\frac{(z-q^2)(z-q^6)(z^2+q^4z+q^8)}{(1-q^2z)(1-q^6z)(1+q^4z+q^8z^2)}(a^{2\La_1})^2.
\end{align*}
\end{proposition}

\section{Fusion construction}

In this section we construct a $U'_q(D_4^{(3)})$-module $V^l$ from $V^1$ by so-called fusion
construction. It is then shown that $V^l$ admits a crystal pseudobase.

\subsection{Review}

Following section 3 of \cite{KMN2} we review the fusion construction and rewrite a necessary 
proposition.

Let $l$ be a positive integer and $\mathfrak{S}_l$ the $l$-th symmetric group.
Let $s_i$ be the simple reflection which interchanges $i$ and $i+1$, and 
let $l(w)$ be the length of $w\in\mathfrak{S}_l$. Let $\geh$ be an affine Lie algebra and $V$
a finite-dimensional $U'_q(\geh)$-module which has a $U'_q(\geh)_{\Kz}$-submodule $V_{\Kz}$.
Assume that $V$ has a polarization $(\;,\;)$ such that $(V_{\Kz},V_{\Kz})\subset\Kz$.
Assume also that $V$ admits a crystal base which is perfect of level $1$.
Let $R(x,y)$ denote the 
$R$-matrix for $V\ot V$. For any $w\in\mathfrak{S}_l$ we construct a 
$U_q'(\geh)$-linear map 
$R_w(x_1,\ldots,x_l):V_{x_1}\ot\cd\ot V_{x_l}\rightarrow V_{x_{w(1)}}\ot\cd\ot V_{x_{w(l)}}$
by
\begin{align*}
R_1(x_1,\ldots,x_l)&=1,\\
R_{s_i}(x_1,\ldots,x_l)&=\left(\bigotimes_{j<i}\text{id}_{V_{x_j}}\right)\ot
R(x_i,x_{i+1})\ot\left(\bigotimes_{j>i+1}\text{id}_{V_{x_j}}\right),\\
R_{ww'}(x_1,\ldots,x_l)&=R_{w'}(x_{w(1)},\ldots,x_{w(l)})\circ R_w(x_1,\ldots,x_l)\\
&\hspace{1cm}\text{ for $w,w'$ such that $l(ww')=l(w)+l(w')$.}
\end{align*}
Fix $r\in\Z_{>0}$. For each $l\in\Z_{>0}$, we put
\begin{align*}
R_l=&R_{w_0}(q^{r(l-1)},q^{r(l-3)},\ldots,q^{-r(l-1)}):\\
&V_{q^{r(l-1)}}\ot V_{q^{r(l-3)}}\ot\cd\ot V_{q^{-r(l-1)}}\rightarrow
V_{q^{-r(l-1)}}\ot V_{q^{-r(l-3)}}\ot\cd\ot V_{q^{r(l-1)}},
\end{align*}
where $w_0$ is the longest element of $\mathfrak{S}_l$. Then $R_l$ is a $U'_q(\geh)$-linear
homomorphism. Define 
\[
V^l=\mbox{Im}\;R_l.
\]
Let us denote by $W$ the image of 
\[
R(q^r,q^{-r}):V_{q^r}\ot V_{q^{-r}}\longrightarrow V_{q^{-r}}\ot V_{q^r}
\]
and by $N$ its kernel. Then we have
\begin{align*}
&V^l\text{ considered as a submodule of }V^{\ot l}=V_{q^{-r(l-1)}}\ot\cd\ot V_{q^{r(l-1)}}\\
&\text{is contained in }\bigcap_{i=0}^{l-2}V^{\ot i}\ot W\ot V^{\ot(l-2-i)}.\\
\intertext{Similarly, we have}
&V^l\text{ is a quotient of }V^{\ot l}/\sum_{i=0}^{l-2}V^{\ot i}\ot N\ot V^{\ot(l-2-i)}.
\end{align*}

Let $P$ be the weight lattice of $\geh$ and set $P_{cl}=P/\Z\delta$. Let $\la_0$ be an 
element of $P_{cl}$ such that
\begin{equation} \label{eq:wt assump}
\wt V\subset\la_0+\sum_{i\neq0}\Z_{\le0}\alpha_i\text{ and }\dim V_{\la_0}=1.
\end{equation}
Take a non-zero vector $u_0$ from $V_{\la_0}$. Let $\vphi(z)$ be a function such that 
\[
R(x,y)(u_0\ot u_0)=\vphi(x/y)(u_0\ot u_0).
\]
We assume that 
\begin{equation} \label{eq:phi assump}
\vphi(q^{2kr})\text{ does not vanish for any }k>0.
\end{equation}
Let $I$ be the index set of the simple roots of $\geh$ and $\geh_{I\setminus\{0\}}$
the finite-dimensional simple Lie algebra whose Dynkin diagram is obtained by removing 
the 0-vertex of that of $\geh$. Let $V(\la)$ be the irreducible 
$U_q(\geh_{I\setminus\{0\}})$-module with highest weight $\la$.

\begin{proposition}{\rm (Proposition 3.4.5 of \cite{KMN2})} \label{prop:3.4.5}
Let $m$ be a positive integer and assume the following conditions:
\begin{itemize}
\item[(i)] $\langle h_i,l\la_0+j\alpha_0\rangle\ge0$ for $i\neq 0$ and $0\le j\le m$.
\item[(ii)] $\dim(V^l)_{l\la_0+k\alpha_0}\le\sum_{j=0}^m
	\dim V(l\la_0+j\alpha_0)_{l\la_0+k\alpha_0}$ for $0\le k\le m$.
\item[(iii)] There exists $i_1\in I$ such that $\{i\in I\mid\langle h_0,\alpha_i\rangle<0\}
	=\{i_1\}$.
\item[(iv)] $-\langle h_0,l\la_0-\alpha_{i_1}\rangle\ge0$.
\end{itemize}
Then we have 
\[
V^l\simeq\bigoplus_{j=0}^m V(l\la_0+j\alpha_0)\quad
\text{as a $U_q(\geh_{I\setminus\{0\}})$-module}
\]
and $V^l$ admits a crystal pseudobase as a $U'_q(\geh)$-module.
\end{proposition}

\subsection{Our case}

Set $\geh=D_4^{(3)}$ and let $V$ be the representation $V^1$ constructed in section 
\ref{subsec:fund rep}. We have checked that $V^1$ has a polarization such that 
$(V^1_{\Kz},V^1_{\Kz})\subset\Kz$. 
We have also calculated the $R$-matrix for $V^1\ot V^1$. Set $r=1$ and $\la_0=\olLa_1$.
Then (\ref{eq:wt assump}) is satisfied. 
{}From Proposition \ref{prop:R mat} we have
\[
\vphi(z)=(1-q^2z)(1-q^6z)(1+q^4z+q^8z^2).
\]
Hence (\ref{eq:phi assump}) is also satisfied.

\begin{theorem} \label{th:V^l}
The $U'_q(D_4^{(3)})$-module $V^l$ constructed by the fusion construction admits a 
crystal pseudobase. Moreover, we have
\[
V^l\simeq\bigoplus_{j=0}^l V(j\olLa_1)\quad\text{as a $U_q(G_2)$-module}.
\]
\end{theorem}

\begin{proof}
We use Proposition \ref{prop:3.4.5}. It suffices to check the conditions (i)-(iv).
Set $\la_0=\olLa_1,m=l$. Note that $\olLa_1=-\alpha_0$. (i), (iii) and (iv) are easily
checked as
\begin{itemize}
\item[(i)] $\langle h_i,l\la_0+j\alpha_0\rangle=(l-j)\delta_{i1}\ge0$ for $i\neq 0$ 
	and $0\le j\le l$.
\item[(iii)] $\{i\in I\mid\langle h_0,\alpha_i\rangle<0\}=\{1\}$.
\item[(iv)] $-\langle h_0,l\la_0-\alpha_1\rangle=2l-1$.
\end{itemize}

We are to show (ii). By the direct calculation using Proposition \ref{prop:R mat}, we 
see $N={\rm Ker}\;R(q,q^{-1})$ 
contains $u_{\La_2},u_{\La_1}^{(1)}-u_{\La_1}^{(2)},
q(1-q^4)u_{\La_1}^{(1)}-u_{\La_1}^{(3)},[2](1-q^4+q^8)u_0^{(1)}-[3]u_0^{(2)}$. Hence
by the explicit form of the highest weight vectors, at $q=1$, $N$ contains $\bigwedge^2
V(\olLa_1),V(0)\wedge V(\olLa_1)$ and $v_\phi\ot v_\phi+u$, where $u$ is an element of
$V(\olLa_1)^{\ot2}$. Hence, at $q=1$, 
\[
V^{\ot l}/\sum V^{\ot j}\ot N\ot V^{\ot(l-2-j)}
\]
is generated by $S^l(V(\olLa_1))$ and $V(0)\ot S^{l-1}(V(\olLa_1))$.
Hence so is for a generic $q$. Thus we obtain
\begin{align*}
\sum_{l\ge0}\mbox{ch}\;(V^l)t^l&
\le\sum_{l\ge0}\mbox{ch}\;S^l(V(\olLa_1))t^l+\sum_{l\ge1}\mbox{ch}\;S^{l-1}(V(\olLa_1))t^l\\
&=(1+t)\sum_{l\ge0}\mbox{ch}\;S^l(V(\olLa_1))t^l\\
&=\frac{1+t}{(1-t)\prod_{\beta\in S}(1-e^\beta t)(1-e^{-\beta}t)},
\end{align*}
where $S=\{\alpha_1,\alpha_1+\alpha_2,2\alpha_1+\alpha_2\}$. On the other hand, from the
description of the crystal base of $U_q(G_2)$-module $V(j\olLa_1)$ in section \ref{subsec:def},
\eq
&&\ba{rl}
\sum_{l\ge0}\left(\sum_{0\le j\le l}\mbox{ch}\;V(j\olLa_1)\right)t^l
&=\dfrac1{1-t}\sum_{j\ge0}\mbox{ch}\;V(j\olLa_1)t^j \\[1ex]
&=\dfrac{1+t}{(1-t)\prod_{\beta\in S}(1-e^\beta t)(1-e^{-\beta}t)}. \ea
\label{eq:gen ch}
\eneq
Thus we have 
\[
\dim(V^l)_\la\le\sum_{j=0}^l\dim V(j\olLa_1)_\la
\]
for any $\la$. The proof is completed.
\end{proof}

\section{$U'_q(D_4^{(3)})$-crystal}

In this section we define a $U'_q(D_4^{(1)})$-crystal $B_l$. 
As a $U_q(G_2)$-crystal, $B_l$ is
isomorphic to the crystal $\bigoplus_{j=0}^lB^{G_2}(j\olLa_1)$ 
for the $U_q(G_2)$-module 
$\bigoplus_{j=0}^lV^{G_2}(j\olLa_1)$.

\subsection{$U_q(G_2)$-crystal} \label{subsec:def}
In \cite{KM} the crystal graph for any finite-dimensional irreducible $U_q(G_2)$-module
was given. For our purpose the description of the crystal $B^{G_2}(j\La_1)$ 
for the highest weight
module with highest weight $j\La_1$ is necessary. Any element of $B^{G_2}(j\La_1)$ is represented
as a one-row semistandard tableau whose entries are 
$1, 2, 3, 0, {\bar 3}, {\bar 2}, {\bar 1}$ with the total order 
$1\prec 2\prec 3\prec 0\prec {\bar 3}\prec {\bar 2}\prec {\bar 1}$ as
\[
\underbrace{
\begin{array}{c}
\hline
\multicolumn{1}{|c|}{\lw{1}\ldots \lw{1}} \\
\hline 
\end{array}
}_{w_1}\!
\underbrace{
\begin{array}{c}
\hline
\multicolumn{1}{|c|}{\lw{2}\ldots \lw{2}} \\
\hline 
\end{array}
}_{w_2}\!
\underbrace{
\begin{array}{c}
\hline
\multicolumn{1}{|c|}{\lw{3}\ldots \lw{3}} \\
\hline 
\end{array}
}_{w_3}\!
\underbrace{
\begin{array}{c}
\hline
\multicolumn{1}{|c|}{\ \lw{0}\ } \\
\hline 
\end{array}
}_{w_0}\!
\underbrace{
\begin{array}{c}
\hline
\multicolumn{1}{|c|}{\lw{\bar 3}\ldots \lw{\bar 3}} \\\hline 
\end{array}
}_{{\bar w}_3}\!
\underbrace{
\begin{array}{c}
\hline
\multicolumn{1}{|c|}{\lw{\bar 2}\ldots \lw{\bar 2}} \\
\hline 
\end{array}
}_{{\bar w}_2}\!
\underbrace{
\begin{array}{c}
\hline 
\multicolumn{1}{|c|}{\lw{\bar 1}\ldots \lw{\bar 1}} \\
\hline 
\end{array}
}_{{\bar w}_1}.
\]
In the tableau,
$0$ occurs at most once and the length is $j$, i.e., $w_0=0$ or $1$,
$\sum_{i=1}^3(w_i+\bar{w}_i)+w_0=j$. For instance 
$
\begin{array}{ccccccc}
\hline
\multicolumn{1}{|c}{\lw{1}} &
\multicolumn{1}{|c}{\lw{2}} &
\multicolumn{1}{|c}{\lw{2}} &
\multicolumn{1}{|c}{\lw{3}} &
\multicolumn{1}{|c}{\lw{0}} &
\multicolumn{1}{|c}{\lw{\bar 1}} &
\multicolumn{1}{|c|}{\lw{\bar 1}} \\
\hline 
\end{array}
$
is an element of $B^{G_2}(7\La_1)$. It is also useful to introduce
a coordinate representation for
an element of $B^{G_2}(j\La_1)$ by 
\begin{align*}
x_i&=w_i,\quad \bar{x}_i=\bar{w}_i\;(i=1,2),\\
x_3&=2w_3+w_0,\quad \bar{x}_3=2\bar{w}_3+w_0.
\end{align*}
Then we have 
\[
B^{G_2}(j\La_1)=\left\{b=(x_1,x_2,x_3,{\bar x}_3,{\bar x}_2,{\bar x}_1)\in(\Zn)^6
\left\vert
\begin{array}{l}
x_3\equiv\bar{x}_3\;(\text{mod }2), \\
\sum_{i=1,2} (x_i+{\bar x}_i)+(x_3+{\bar x}_3)/2=j
\end{array}
\right.
\right\}. 
\]
Below we give the explicit crystal structure of $B^{G_2}(j\La_1)$ with this parametrization. Set 
$(x)_+=\max(x,0)$, then we have 
\begin{align*}
{\tilde e}_1 b=&
\begin{cases}
(\ldots,{\bar x}_2 +1,{\bar x}_1 -1) 
& \text{if ${\bar x}_2 -{\bar x}_3 \geq (x_2 -x_3)_+$}, 
\\  
(\ldots,x_3 +1,{\bar x}_3 -1,\ldots) 
& \text{if ${\bar x}_2 -{\bar x}_3 <0\leq x_3 -x_2$}, 
\\ 
(x_1 +1,x_2 -1,\ldots) 
& \text{if $({\bar x}_2 -{\bar x}_3)_+ <x_2 -x_3$},
\end{cases}
\\
{\tilde f}_1 b=&
\begin{cases}
(x_1 -1,x_2 +1,\ldots) 
& \text{if $({\bar x}_2 -{\bar x}_3)_+ \leq x_2 -x_3$}, 
\\
(\ldots,x_3 -1,{\bar x}_3 +1,\ldots) 
& \text{if ${\bar x}_2 -{\bar x}_3 \leq 0<x_3 -x_2$}, 
\\
(\ldots,{\bar x}_2 -1,{\bar x}_1 +1) 
& \text {if ${\bar x}_2 -{\bar x}_3 >(x_2 -x_3)_+$},
\end{cases}
\\
{\tilde e}_2 b=&
\begin{cases}
(\ldots,{\bar x}_3 +2,{\bar x}_2 -1,\ldots) 
& \text{if ${\bar x}_3 \geq x_3$}, 
\\
(\ldots,x_2 +1,x_3 -2,\ldots) 
& \text{if ${\bar x}_3 <x_3$},
\end{cases}
\\
{\tilde f}_2 b=&
\begin{cases}
(\ldots,x_2 -1,x_3 +2,\ldots) 
& \text{if ${\bar x}_3 \leq x_3$}, 
\\
(\ldots,{\bar x}_3 -2,{\bar x}_2 +1,\ldots) 
& \text{if ${\bar x}_3 >x_3$},
\end{cases}
\end{align*}
\begin{align*}
\veps_1(b)=&{\bar x}_1+({\bar x}_3-{\bar x}_2+(x_2-x_3)_+)_+,
&
\veps_2(b)=&{\bar x}_2+\frac{1}{2}(x_3-{\bar x}_3)_+,
\\
\varphi_1(b)=&x_1+(x_3-x_2+({\bar x}_2-{\bar x}_3)_+)_+,
&
\varphi_2(b)=&x_2+\frac{1}{2}({\bar x}_3-x_3)_+.
\end{align*}
If $\et{i}b$ or $\ft{i}b$ does not belong to $B^{G_2}(j\La_1)$, namely, if $x_j$ or $\bar{x}_j$
for some $j$ becomes negative, we should understand it to be $0$.

\subsection{Action of $\et{0},\ft{0}$} \label{subsec:0-action}

For a positive integer $l$ we introduce a 
$U_q'(\geh)$-crystal $B_l$.
As a $U_q(G_2)$-crystal,
\[
B_l=\bigoplus_{j=0}^lB^{G_2}(j\La_1),
\]
where $B^{G_2}(j\La_1)$ is the $U_q(G_2)$-crystal 
explained in the previous subsection. 
To define the actions of $\et{0}$ and $\ft{0}$, we
introduce conditions ($E_1$)-($E_6$) and ($F_1$)-($F_6$). Set
\begin{equation} \label{z1-4}
z_1={\bar x}_1-x_1, \quad 
z_2={\bar x}_2 -{\bar x}_3, \quad 
z_3=x_3-x_2, \quad 
z_4=({\bar x}_3-x_3)/2,
\end{equation}
and
\begin{eqnarray*}
&&(F_1)\quad
z_1+z_2+z_3+3z_4\le0, z_1+z_2+3z_4\le0, z_1+z_2\le0, z_1\le0,
\\
&&(F_2)\quad
z_1+z_2+z_3+3z_4\le0, z_2+3z_4\le0, z_2\le0, z_1> 0,
\\
&&(F_3)\quad 
z_1+z_3+3z_4\le0, z_3+3z_4\le0, z_4\le0, z_2> 0, z_1+z_2> 0,
\\
&&(F_4)\quad
z_1+z_2+3z_4> 0, z_2+3z_4> 0, z_4> 0, z_3\le0, z_1+z_3\le0,
\\
&&(F_5)\quad
z_1+z_2+z_3+3z_4> 0, z_3+3z_4> 0, z_3> 0, z_1\le0, 
\\
&&(F_6)\quad
z_1+z_2+z_3+3z_4> 0, z_1+z_3+3z_4> 0, z_1+z_3> 0, z_1> 0.
\end{eqnarray*}
($E_i$) ($1\le i\le 6$) is defined from ($F_i$) by replacing $>$ (resp. $\le$) with 
$\ge$ (resp. $<$). We define
\begin{align*}
{\tilde e}_0 b=&
\begin{cases}
\mathscr{E}_1 b:=
(x_1 -1,\ldots) 
& \text{if ($E_1$)}, 
\\
\mathscr{E}_2 b:=
(\ldots,x_3 -1,{\bar x}_3 -1,\ldots,{\bar x}_1 +1) 
& \text{if ($E_2$)}, 
\\
\mathscr{E}_3 b:=
(\ldots,x_3 -2,\ldots,{\bar x}_2 +1,\ldots) 
& \text{if ($E_3$)}, 
\\ 
\mathscr{E}_4 b:=
(\ldots,x_2 -1,\ldots,{\bar x}_3 +2,\ldots) 
& \text{if ($E_4$)}, 
\\
\mathscr{E}_5 b:=
(x_1 -1,\ldots,x_3 +1,{\bar x}_3 +1,\ldots) 
& \text{if ($E_5$)}, 
\\
\mathscr{E}_6 b:=
(\ldots,{\bar x}_1 +1) 
& \text{if ($E_6$)},
\end{cases}
\\
{\tilde f}_0 b=&
\begin{cases}
\mathscr{F}_1 b:=(x_1 +1,\ldots) 
& \text{if ($F_1$)}, 
\\ 
\mathscr{F}_2 b:=
(\ldots,x_3 +1,{\bar x}_3 +1,\ldots,{\bar x}_1 -1) 
& \text{if ($F_2$)}, 
\\
\mathscr{F}_3 b:=
(\ldots,x_3 +2,\ldots,{\bar x}_2 -1,\ldots) 
& \text{if ($F_3$)}, 
\\
\mathscr{F}_4 b:=
(\ldots,x_2 +1,\ldots,{\bar x}_3 -2,\ldots) 
& \text{if ($F_4$)}, 
\\
\mathscr{F}_5 b:=
(x_1 +1,\ldots,x_3 -1,{\bar x}_3 -1,\ldots) 
& \text{if ($F_5$)}, 
\\
\mathscr{F}_6 b:=(\ldots,{\bar x}_1 -1)
& \text{if ($F_6$)}.
\end{cases}
\end{align*}
\begin{remark}
\be[(i)]
\item
Set
\begin{equation} \label{A}
A=(0,z_1,z_1+z_2,z_1+z_2+3z_4,z_1+z_2+z_3+3z_4,2z_1+z_2+z_3+3z_4)
\end{equation}
and $z_1,z_2,z_3,z_4$ are given in \eqref{z1-4}.
Denote the $i$-th component of $A$ by $A_i$. 
Then, for $1\le i\le 6$, $(F_i)$ holds if
and only if $\max A=A_i$ and $A_j<A_i$ for any $j$ such that $1\le j<i$.
Similarly, 
$(E_i)$ holds if
and only if $\max A=A_i$ and $A_j<A_i$ for any $j$ such that $j>i$.
\item
By (i), we have
\[
B_l=
\bigsqcup_{1\leq i\leq 6} \{\text{$b\in B_l$ $\vert$ $b$ satisfies ($E_i$)}\}=
\bigsqcup_{1\leq i\leq 6} \{\text{$b\in B_l$ $\vert$ $b$ satisfies ($F_i$)}\}. 
\]
\ee
\end{remark}

For $b=(x_1,x_2,x_3,\bar{x}_3,\bar{x}_2,\bar{x}_1)$ we set 
\begin{equation} \label{def s(b)}
s(b)=x_1+x_2+\frac{x_3+{\bar x}_3}{2}+{\bar x}_2+{\bar x}_1. 
\end{equation}
Suppose $b\in B_l$. Looking at the rule of the action of $\ft{0}$ carefully, we see that 
the coordinates of $\ft{0}b$ never get negative. 
It means that $\ft{0}b=0$ occurs only when 
$s(b)=l$ and $b$ satisfies $(F_1)$. The case of $\et{0}$ is similar.
Checking directly one can show that $B_l$ satisfies the condition: for $b,b'\in B_l$,
\[
b'=\ft{0}b\Longleftrightarrow b=\et{0}b'.
\]
Hence one can draw the crystal graph of $B_l$ with arrows of color $0,1,2$.

\begin{example}
Let us denote the elements of $B_1$ by 
\begin{align*}
\fbox{$1$}&=(1,0,0,0,0,0),& 
\fbox{$2$}&=(0,1,0,0,0,0),& 
\fbox{$3$}&=(0,0,2,0,0,0),\\
\fbox{$0$}&=(0,0,1,1,0,0),&
\fbox{${\bar 3}$}&=(0,0,0,2,0,0),&    
\fbox{${\bar 2}$}&=(0,0,0,0,1,0),\\
\fbox{${\bar 1}$}&=(0,0,0,0,0,1), &
\phi&=(0,0,0,0,0,0).
\end{align*}
then, the crystal graph of $B_1$ is given as follows: 

\begin{center}
\begin{picture}(220,80)
\put(0,30)
{\makebox(0,0){$\fbox{$1$}\xrightarrow{1}$}}
\put(30,30)
{\makebox(0,0){$\fbox{$2$}\xrightarrow{2}$}}
\put(60,30)
{\makebox(0,0){$\fbox{$3$}\xrightarrow{1}$}}
\put(90,30)
{\makebox(0,0){$\fbox{$0$}\xrightarrow{1}$}}
\put(120,30)
{\makebox(0,0){$\fbox{${\bar 3}$}\xrightarrow{2}$}}
\put(155,30)
{\makebox(0,0){$\fbox{${\bar 2}$}\xrightarrow{1}\fbox{${\bar 1}$}$}}
\put(82,3)
{\makebox(0,0){$\phi$}}
\bezier{200}(86,3)(125,5)(162,20)
\bezier{200}(1,20)(30,5)(77,3)
\bezier{200}(59,38)(100,70)(135,40)
\bezier{200}(30,38)(65,70)(107,40)
\bezier{200}(86,3)(87,3)(88,4)
\bezier{200}(86,3)(88,3)(89,2)
\bezier{200}(1,20)(2,19)(2,18)
\bezier{200}(1,20)(1,20)(3,20)
\bezier{200}(30,38)(31,39)(31,40)
\bezier{200}(30,38)(31,39)(32,39)
\bezier{200}(59,38)(60,39)(60,40)
\bezier{200}(59,38)(60,39)(62,39)
\end{picture}
\end{center}
The arrows without number are $0$-arrows. 
\end{example}

The next two propositions are related to the action of $\et{0},\ft{0}$. The first one is easily 
proved.

\begin{lemma} \label{lem:F1E6}
\begin{itemize}
\item[(1)] Suppose that $b\in B_l$ satisfies $(F_1)$ and $\ft{0}b\in B_l$. Then $\ft{0}b$ also satisfies
	$(F_1)$.
\item[(2)] Suppose that $b\in B_l$ satisfies $(E_6)$ and $\et{0}b\in B_l$. Then $\et{0}b$ also satisfies
	$(E_6)$.
\end{itemize}
\end{lemma}

\begin{proposition} \label{prop:eps0-phi0}
The values of $\veps_0$ and $\vphi_0$ of an element $b=(x_1,x_2,x_3,\bar{x}_3,\bar{x}_2,\bar{x}_1)$
of $B_l$ is given by 
\begin{align*}
\vphi_0(b)&=l-s(b)+\max A, \\
\veps_0(b)&=l-s(b)+\max A-(2z_1+z_2+z_3+3z_4),
\end{align*}
where $A$ is as in \eqref{A}.
\end{proposition}

\begin{proof}
Notice that if $\ft{0}b=0$ 
occurs for $b\in B_l$, then $b$ satisfies $(F_1)$. 
{}From Lemma \ref{lem:F1E6} (1), one verifies that
the formula of $\vphi_0$ is correct when $b$ satisfies $(F_1)$. Thus we are left to show that 
$\vphi_0(\ft{0}b)=\vphi_0(b)-1$ if $b,\ft{0}b\in B_l$. It can be checked case by case. Let $A'$ be 
the list $A$ for $\ft{0}b$ and $A'_i$ be its $i$-th component. Notice that if $b$ satisfies $(F_i)$,
then $\max A'=A'_i$.
\end{proof}

\section{Decomposition of $B_l$ as a $U_q(A_2)$-crystal}

\subsection{Review on $U_q(A_2)$-crystal} \label{subsec:A2 review}

We review on the $U_q(A_2)$-crystal $B^{A_2}(j_0\La_0+j_1\La_1)$ of the highest weight module of
highest weight $j_0\La_0+j_1\La_1$. We use $\{0,1\}$ as the index set of simple roots of $A_2$. 
It is known that any element of $B^{A_2}(j_0\La_0+j_1\La_1)$ is uniquely expressed as 
$\ft{0}^r\ft{1}^q\ft{0}^pu$ for some $p,q,r$ such that 
$0\le p\le j_0$, $p\le q\le j_1+p$, $0\le r\le
j_0-2p+q$. Here $u$ stands for the highest weight vector. By \cite{KN} an element of 
$B^{A_2}(j_0\La_0+j_1\La_1)$ is also represented by a two-row tableau. $\ft{0}^r\ft{1}^q\ft{0}^pu$
corresponds to 
\[
t(p,q,r)=
\begin{array}{l}
1^{j_0+j_1-p-r}\ 2^r\ 3^p \\
2^{j_1+p-q}\ 3^{q-p}
\end{array}.
\]
Below we give the crystal structure.
\begin{align}
{\tilde e}_0t(p,q,r)=&
\begin{cases}
t(p,q,r-1) & \text{if $r>0$}, \\
0          & \text{if $r=0$},
\end{cases}
\label{A2 e0}\\
{\tilde e}_1t(p,q,r)=&
\begin{cases}
t(p-1,q-1,r+1) & \text{if $p>0, p-q+r\geq 0$}, \\
t(p,q-1,r)     & \text{if $p-q+r<0$}, \\
0              & \text{if $p=0, p-q+r\geq 0$},
\end{cases}
\label{A2 e1}\\
{\tilde f}_0t(p,q,r)=&
\begin{cases}
t(p,q,r+1) & \text{if $0\leq r<j_0+q-2p$}, \\
0          & \text{if $r=j_0+q-2p$},
\end{cases}
\label{A2 f0}\\
{\tilde f}_1t(p,q,r)=&
\begin{cases}
t(p+1,q+1,r-1) & \text{if $p\leq q<p+r$},   \\
t(p,q+1,r)     & \text{if $p+r\leq q<j_1+p$}, \\
0              & \text{if $p+r\leq q=j_1+p$}.
\end{cases}
\label{A2 f1}
\end{align}
The remaining data $\veps_i,\vphi_i$ of $t(p,q,r)$ are given by 
\eq
&&\ba{rlrl}\veps_0=&r,\quad &\vphi_0=&j_0-2p+q-r,\\[1ex]
\veps_1=&p+(q-p-r)_+,\quad &\vphi_1=&(p-q+r)_++j_1+p-q. \ea
\label{A2 eps phi}
\eneq

The following proposition is immediate.

\begin{proposition} \label{prop:up side down}
The lowest weight vector of $B^{A_2}(j_0\La_0+j_1\La_1)$ is given by $t(j_0,j_0+j_1,j_1)$.
Moreover, we have 
\[
t(p,q,r)=\et{0}^{r'}\et{1}^{q'}\et{0}^{p'}t(j_0,j_0+j_1,j_1),
\]
where $p'=j_1-q+p,q'=j_0+j_1-q,r'=j_0+q-2p-r$.
\end{proposition}


\subsection{$U_q(A_2)$-crystal structure}

In what follows in this section we investigate the structure of the crystal subgraph of $B_l$ 
obtained by forgetting $2$-arrows.

\begin{define}
For $l\in \mathbb{Z}_{>0}$ take integers $i,j_0,j_1$ such that 
\begin{equation*}
0\leq i\leq {l}/{2},\quad
i\leq j_0,j_1 \leq l-i
\quad\text{and}\quad
j_0,j_1\equiv l-i\kern-1ex\pmod 3,
\end{equation*}
and set $y_a=(l-i-j_a)/3$ for $a=0,1$. We define the element $\bar{b}_{j_0,j_1}^{l,i}$ of $B_l$ by
\begin{equation} \label{def of b-bar}
{\bar b}_{j_0,j_1}^{l,i}=
\begin{cases}
(0,y_1,-2y_1+3y_0+i,y_0+i,y_0+j_0,0)&\text{ if }j_0\le j_1,\\
(0,y_0,y_0+i,2y_1-y_0+i,-y_1+2y_0+j_0,0)&\text{ if }j_0>j_1.
\end{cases}
\end{equation}
We also define the subset $B_{j_0,j_1}^{l,i}$ of $B_l$ to be the connected component of $B_l$ 
generated by $\et{a},\ft{a}$ $(a=0,1)$ that contains $\bar{b}_{j_0,j_1}^{l,i}$.
\end{define}

Our main theorem of this section is given as follows.

\begin{theorem}\label{thm:$U_q(A_2)$CrystalStructure}
Forgetting $2$-arrows, the crystal graph $B_l$ decomposes into connected components in the 
following manner.
\[
B_l= 
{\bigsqcup}^{[\frac{l}{2}]}_{i=0} 
{\bigsqcup}_{\begin{subarray}{c}
{i\leq j_0, j_1 \leq l-i}\\
{j_0, j_1 \equiv l-i  \!\!\!\!\pmod 3}
\end{subarray}}
B^{l,i}_{j_0,j_1}.
\]
Moreover, $B_{j_0,j_1}^{l,i}$ is isomorphic to the $U_q(A_2)$-crystal $B^{A_2}(j_0\La_0+j_1\La_1)$. 
\end{theorem}

For the proof we introduce some notations. Set 
\[
B_{\ge0}=\left\{(x_1,x_2,x_3,\bar{x}_3,\bar{x}_2,\bar{x}_1)
\in \mathbb{Z}_{\geq 0}^6\ \vert \ 
(x_3+\bar{x}_3)/2\in \mathbb{Z}_{\geq 0} \right\}.
\]
Note that $B_l=\{b\in B_{\ge0}\mid s(b)\le l\}$ where $s(b)$ was defined in \eqref{def s(b)}.
One can endow $B_{\ge0}$ with the crystal structure by applying the same rule for $\et{i},\ft{i}$
as section \ref{subsec:def} and \ref{subsec:0-action} with $l=\infty$. Namely, $\et{i},\ft{i}$ 
vanish only when some coordinate becomes negative. Note that ${\bar b}_{j_0,j_1}^{l,i}$ is 
$U_q(A_2)$-highest, i.e., $\et{a}{\bar b}_{j_0,j_1}^{l,i}=0$ for $a=0,1$, as an element of $B_l$, 
but $\et{0}{\bar b}_{j_0,j_1}^{l,i}\neq0$ as an element of $B_{\ge0}$.

\subsection{Some relations on $B_{\ge0}$}

We prepare two relations that hold on $B_{\ge0}$.

\begin{lemma} \label{lem:onion}
Suppose that $j_0\le j_1$. On $B_{\ge0}$ we have 
\[
\ft{0}\ft{1}^q\ft{0}^p{\bar b}_{j_0,j_1}^{l,i}=\ft{1}^{q-1}\ft{0}^p{\bar b}_{j_0-1,j_1-1}^{l-1,i}
\mbox{ if $i<j_0$, $p<j_0$, $p<q\le j_1+p$.}
\]
\end{lemma}

\begin{proof}
We use the table in Appendix \ref{app:A}. Under the assumption one can show that all the cases
satisfy $(F_6)$ of the rule of 0-action. Hence, if we write 
$x=\ft{1}^q\ft{0}^p{\bar b}_{j_0,j_1}^{l,i}=(x_1,x_2,x_3,{\bar x}_3,{\bar x}_2,{\bar x}_1)$,
then $\ft{0}x=(x_1,x_2,x_3,{\bar x}_3,{\bar x}_2,{\bar x}_1-1)$. On the other hand, in each case
we also have
\[
x\left|_{(j_0,j_1,q)\rightarrow(j_0-1,j_1-1,q-1)}\right.
=(x_1,x_2,x_3,{\bar x}_3,{\bar x}_2,{\bar x}_1-1).
\]
Hence we have the desired relation. Note that $l$ should be replaced by $l-1$ so that $y_0,y_1$
remain the same.
\end{proof}

\begin{lemma} \label{lem:comm}
Suppose that $j_0\le j_1$. On $B_{\ge0}$ we have 
\[
\ft{0}\ft{1}^q\ft{0}^p{\bar b}_{j_0,j_1}^{l,i}=\ft{1}^q\ft{0}^{p+1}{\bar b}_{j_0,j_1}^{l,i}
\mbox{ if }q\le p<j_0.
\]
\end{lemma}

\begin{proof}
We again use the table in Appendix \ref{app:A}. Under the assumption the cases that occur are
\begin{itemize}
\item[(1)] $I.$ $0\le p\le i$, (i) $0\le q\le j_0-i+p$,
\item[(2)] $II.$ $i\le p$, (i) $0\le q\le j_0-p+i$,
\item[(3)] $II.$ $i\le p$, (ii) $j_0-p+i\le q\le j_1+p-i$.
\end{itemize}
Each case satisfies $(F_5),(F_3),(F_2)$, respectively. In each case the action of $\ft{0}$ is realized 
by replacing $p$ with $p+1$.
\end{proof}

\subsection{Proof of Theorem \ref{thm:$U_q(A_2)$CrystalStructure}}

We frequently use the following condition for $(p,q,r)$.
\[
(C)\quad 0\le p\le j_0,\ p\le q\le j_1+p,\ 0\le r\le j_0+q-2p.
\]

\noindent[Step 1] We show
\[
B_{j_0,j_1}^{l,i}\simeq B^{A_2}(j_0\La_0+j_1\La_1)\text{ for }j_0\le j_1.
\]

\begin{proof}
Due to the fact that for $b,b'\in B_l$, $b'=\ft{i}b$ if and only if $b=\et{i}b'$ ($i=0,1$), 
it suffices to show for $(C)$
\begin{itemize}
\item[(i)] $\ft{0}^r\ft{1}^q\ft{0}^p\bar{b}_{j_0,j_1}^{l,i}\in B_l$,

\item[(ii)] $b=\ft{0}^{j_0+q-2p}\ft{1}^q\ft{0}^p\bar{b}_{j_0,j_1}^{l,i}$ satisfies $(F_1)$ and 
$s(b)=l$,
\end{itemize}
and as an element of $B_l$
\begin{itemize}
\item[(iii)] $\tilde{f}_1\tilde{f}_0^r\tilde{f}_1^q\tilde{f}_0^p \bar{b}_{j_0,j_1}^{l,i}=
\begin{cases}
\tilde{f}_0^{r-1}\tilde{f}_1^{q+1}\tilde{f}_0^{p+1} \bar{b}_{j_0,j_1}^{l,i} & \text{if $p\leq q<p+r$},
\\ 
\tilde{f}_0^r\tilde{f}_1^{q+1}\tilde{f}_0^p \bar{b}_{j_0,j_1}^{l,i} & \text{if $p+r\leq q<j_1+p$},
\\
0 & \text{if $p+r\leq q=j_1+p$},
\end{cases}$

\item[(iv)] $\et{0}\ft{1}^q\ft{0}^p \bar{b}_{j_0,j_1}^{l,i}=0$,

\item[(v)] $\et{1}\ft{0}^r\ft{1}^q \bar{b}_{j_0,j_1}^{l,i}=0\text{ if }r\ge q$.
\end{itemize}
We prove (i)-(v) by using induction on $l$.

(i) Suppose that $r>0,i<j_0,p<j_0,p<q$. By Lemma \ref{lem:onion} we have 
\[
\ft{0}^r\ft{1}^q\ft{0}^p\bar{b}_{j_0,j_1}^{l,i}
=\ft{0}^{r-1}\ft{1}^{q-1}\ft{0}^p\bar{b}_{j_0-1,j_1-1}^{l-1,i}.
\]
If $r<j_0+q-2p$, we get $\ft{0}^{r-1}\ft{1}^{q-1}\ft{0}^p\bar{b}_{j_0-1,j_1-1}^{l-1,i}\in B_{l-1}
\subset B_l$ by induction hypothesis. If $r=j_0+q-2p$, we know that 
\[
\ft{0}^{j_0+q-2p-2}\ft{1}^{q-1}\ft{0}^p\bar{b}_{j_0-1,j_1-1}^{l-1,i}\mbox{ satisfies $(F_1)$ and }
s(b)=l-1
\]
by induction hypothesis. Hence $\ft{0}^{j_0+q-2p-1}\ft{1}^{q-1}\ft{0}^p
\bar{b}_{j_0-1,j_1-1}^{l-1,i}\in B_l$.

In the cases of $r=0,i=j_0,p=j_0$, $\ft{0}^r\ft{1}^q\ft{0}^p\bar{b}_{j_0,j_1}^{l,i}\in B_l$ can
be checked directly by consulting the table in Appendix \ref{app:A},\ref{app:B},\ref{app:C},
respectively. If $p=q$, we have 
\[
\ft{0}^r\ft{1}^p\ft{0}^p\bar{b}_{j_0,j_1}^{l,i}=\ft{1}^p\ft{0}^{p+r}\bar{b}_{j_0,j_1}^{l,i}
\]
by Lemma \ref{lem:comm} and $\ft{1}^p\ft{0}^{p+r}\bar{b}_{j_0,j_1}^{l,i}\in B_l$ can be checked
directly from Appendix \ref{app:A}.

(ii) The claim can be checked directly from Appendix \ref{app:D}.

(iii) Suppose that $r>0$, $i<j_0$, $p<j_0,p<q$. 
By Lemma \ref{lem:onion} we have 
\[
\ft{1}\ft{0}^r\ft{1}^q\ft{0}^p\bar{b}_{j_0,j_1}^{l,i}
=\ft{1}\ft{0}^{r-1}\ft{1}^{q-1}\ft{0}^p\bar{b}_{j_0-1,j_1-1}^{l-1,i}.
\]
If $r<j_0+q-2p$, the claim is proved by using induction hypothesis and Lemma \ref{lem:onion} once 
again. If $r=j_0+q-2p$, one needs to show 
$\ft{1}\ft{0}^{j_0+q-2p}\ft{1}^q\ft{0}^p\bar{b}_{j_0,j_1}^{l,i}
=\ft{0}^{j_0+q-2p-1}\ft{1}^{q+1}\ft{0}^{p+1}\bar{b}_{j_0,j_1}^{l,i}$, which can be checked from
Appendix \ref{app:D}.

If $r=0$, one needs to show
\[
\ft{1}\ft{1}^q\ft{0}^p \bar{b}_{j_0,j_1}^{l,i}=
\begin{cases}
\ft{1}^{q+1}\ft{0}^p \bar{b}_{j_0,j_1}^{l,i} & \text{if $q<j_1+p$},\\
0 & \text{if $q=j_1+p$},
\end{cases}
\]
which can be checked from Appendix \ref{app:A}. If $i=j_0$, the claim can be checked from Appendix
\ref{app:B}. If $p=j_0$, one needs to show 
\[
\ft{1}\ft{0}^r\ft{1}^q\ft{0}^{j_0} \bar{b}_{j_0,j_1}^{l,i}=
\begin{cases}
\ft{0}^r\ft{1}^{q+1}\ft{0}^{j_0} \bar{b}_{j_0,j_1}^{l,i} & \text{if $q<j_1+j_0$},\\
0 & \text{if $q=j_1+j_0$},
\end{cases}
\]
which can be checked from Appendix \ref{app:C}. If $p=q$ and $r>0$, we have
\[
\ft{1}\ft{0}^r\ft{1}^p\ft{0}^p \bar{b}_{j_0,j_1}^{l,i}=
\ft{1}^{p+1}\ft{0}^{p+r}\bar{b}_{j_0,j_1}^{l,i}=
\ft{0}^{r-1}\ft{1}^{p+1}\ft{0}^{p+1} \bar{b}_{j_0,j_1}^{l,i}
\]
by Lemma \ref{lem:comm}.

(iv) The claim is checked directly from Appendix \ref{app:A}.

(v) Suppose that $r>0$, $j_0>0$, $i<j_0,q>0$. By Lemma \ref{lem:onion} we have 
\[
\et{1}\ft{0}^r\ft{1}^q\bar{b}_{j_0,j_1}^{l,i}=
\et{1}\ft{0}^{r-1}\ft{1}^{q-1}\bar{b}_{j_0-1,j_1-1}^{l-1,i}.
\]
If $r<j_0+q$, the RHS is $0$ by induction hypothesis. If $r=j_0+q$, 
$\et{1}\ft{0}^{j_0+q}\ft{1}^q\bar{b}_{j_0,j_1}^{l,i}=0$ can be checked directly from Appendix
\ref{app:D}.

Note that $r=0$ implies $q=0$ and $j_0=0$ implies $i=j_0$. If $i=j_0$, the claim is checked from 
Appendix \ref{app:B}. If $q=0$, it is checked from Appendix \ref{app:A}.
\end{proof}

\noindent[Step 2] 
Next we show 
\[
B_{j_0,j_1}^{l,i}\simeq B^{A_2}(j_0\La_0+j_1\La_1)\text{ for }j_0>j_1.
\]

Define an involution on $B_{\ge0}$ by
\[
b=(x_1,x_2,x_3,\bar{x}_3,\bar{x}_2,\bar{x}_1) \ \mapsto \ 
(\bar{x}_1,\bar{x}_2,\bar{x}_3,x_3,x_2,x_1)=b^{\vee}. 
\]
We prove two lemmas related to this involution. The next one follows immediately from
the definitions.

\begin{lemma} \label{lem:invol}
Let $b\in B_l$. For $i=0,1,2$, 
\begin{itemize}
\item[(1)] if $\et{i}b\neq0$, then $(\et{i}b)^\vee=\ft{i}(b^\vee)$.
\item[(2)] if $\ft{i}b\neq0$, then $(\ft{i}b)^\vee=\et{i}(b^\vee)$.
\end{itemize}
\end{lemma}

\begin{lemma} \label{lem:invol2}
Suppose that $j_0\le j_1$. As an element of $B_l$, we have for (C)
\[
(\ft{0}^r\ft{1}^q\ft{0}^p \bar{b}_{j_0,j_1}^{l,i})^\vee
=\ft{0}^{r'}\ft{1}^{q'}\ft{0}^{p'} \bar{b}_{j_1,j_0}^{l,i},
\]
where
\begin{equation} \label{p'q'r'}
p'=j_1-q+p,\;q'=j_0+j_1-q,\;r'=j_0+q-2p-r.
\end{equation}
\end{lemma}

\begin{proof}
By the result of Step 1 and Proposition \ref{prop:up side down} we have 
$\ft{0}^r\ft{1}^q\ft{0}^p \bar{b}_{j_0,j_1}^{l,i}=
\et{0}^{r'}\et{1}^{q'}\et{0}^{p'}{\underline b}_{j_0,j_1}^{l,i}$, where 
${\underline b}_{j_0,j_1}^{l,i}=\ft{0}^{j_1}\ft{1}^{j_0+j_1}\ft{0}^{j_0}\bar{b}_{j_0,j_1}^{l,i}$.
Apply ${}^\vee$ on both sides and use the previous lemma. We obtain  
$({\tilde f}_0^r{\tilde f}_1^q{\tilde f}_0^p{\bar b}_{j_0,j_1}^{l,i})^\vee=
\ft{0}^{r'}\ft{1}^{q'}\ft{0}^{p'}({\underline b}_{j_0,j_1}^{l,i})^\vee$, which can be shown from 
the table in Appendix \ref{app:D} and the definition \eqref{def of b-bar}.
\end{proof}

\noindent{\it Proof of Step 2.}
Apply ${}^\vee$ on both sides of (i)-(v) in the proof of Step 1. Use 
Lemma \ref{lem:invol2} and interchange $(p,q,r)$ and $(p',q',r')$. Substituting $(p',q',r')$ with 
$j_0$ and $j_1$ interchanged we obtain for $0\le p\le j_1,p\le q\le j_0+p,0\le r\le j_1-2p+q$,

\begin{itemize}
\item[(i')] $\ft{0}^r\ft{1}^q\ft{0}^p\bar{b}_{j_1,j_0}^{l,i}\in B_l$,

\item[(ii')] $\et{0}\tilde{f}_0^r\tilde{f}_1^q\tilde{f}_0^p \bar{b}_{j_1,j_0}^{l,i}=
\begin{cases}
\tilde{f}_0^{r-1}\tilde{f}_1^q\tilde{f}_0^p \bar{b}_{j_1,j_0}^{l,i} & \text{if $r>0$},
\\
0 & \text{if $r=0$},
\end{cases}$

\item[(iii')] $\et{1}\tilde{f}_0^r\tilde{f}_1^q\tilde{f}_0^p \bar{b}_{j_1,j_0}^{l,i}=
\begin{cases}
\tilde{f}_0^{r}\tilde{f}_1^{q-1}\tilde{f}_0^{p} \bar{b}_{j_1,j_0}^{l,i} & \text{if $p-q+r<0$},
\\ 
\tilde{f}_0^{r+1}\tilde{f}_1^{q-1}\tilde{f}_0^{p-1} \bar{b}_{j_1,j_0}^{l,i} 
& \text{if $p>0,p-q+r\ge0$},
\\
0 & \text{if $p=0,p-q+r\ge0$},
\end{cases}$

\item[(iv')] $\ft{0}\ft{0}^{j_1+q-2p}\ft{1}^q\ft{0}^p \bar{b}_{j_1,j_0}^{l,i}=0$,

\item[(v')] $\ft{1}\ft{0}^r\ft{1}^q \bar{b}_{j_1,j_0}^{l,i}=0\text{ if }p+r\le q=j_0+p$.
\end{itemize}
These relations are enough to check our claim.
\qed

\medskip
\noindent[Step 3] 
We are left to show
\[
B_l= 
{\bigsqcup}^{[\frac{l}{2}]}_{i=0} 
{\bigsqcup}_{\begin{subarray}{c}
{i\leq j_0, j_1 \leq l-i}\\
{j_0, j_1 \equiv l-i  \!\!\!\!\pmod 3}
\end{subarray}}
B^{l,i}_{j_0,j_1}. 
\]

\begin{proof} 
It suffices to show 
\[
{\sharp B}_l={\sum}^{[\frac{l}{2}]}_{i=0} 
{\sum}_{\begin{subarray}{c}
{i\leq j_0, j_1 \leq l-i}\\
{j_0, j_1 \equiv l-i  \!\!\!\!\pmod 3}
\end{subarray}}
{\sharp B}^{l,i}_{j_0,j_1}. 
\]
By Step 1 and 2 we have 
\[
{\sharp B}_{j_0,j_1}^{l,i}=
{\sharp B}^{A_2}(j_0\La_0+j_1\La_1)=
\frac{(1+j_0)(1+j_1)(2+j_0 +j_1)}{2}. 
\]
By direct calculation we have 
\[
{\sum}^{[\frac{l}{2}]}_{i=0} 
{\sum}_{\begin{subarray}{c}
{i\leq j_0, j_1 \leq l-i}\\
{j_0, j_1 \equiv l-i  \!\!\!\!\pmod 3}
\end{subarray}}
{\sharp B}^{l,i}_{j_0,j_1}
=\frac{(l+1)(l+2)(l+3)^2(l+4)(l+5)}{360}. 
\]
On the other hand, computing ${\sharp B}_l$ from the definition of $B_l$ reads
\[
{\sharp B_l}=\frac{(l+6)!}{l!6!}+\frac{((l-1)+6)!}{(l-1)!6!}=\frac{(l+1)(l+2)(l+3)^2(l+4)(l+5)}{360}.
\]
\end{proof}

Thus the proof of Theorem \ref{thm:$U_q(A_2)$CrystalStructure} is completed.

\section{Uniqueness problem}

In this section we deal with a certain uniqueness problem 
of crystals in a more general situation.
In order to state our theorem precisely, 
we prepare some notations. Let $\geh$ be an affine Lie
algebra and $I$ the index set of vertices of the corresponding Dynkin diagram. Let 
$\{\alpha_i\}_{i\in I},\{h_i\}_{i\in I},\{\La_i\}_{i\in I}$ be the set of simple roots, simple
coroots, fundamental weights. Let $\delta$ and $c$ be the generator of 
null roots and the 
canonical central element. 
Let $0$ be the vertex of the Dynkin diagram as in \cite{Kac}.
Let $(a_{ij})_{i,j\in I}$ be the generalized Cartan matrix of $\geh$.
We assume the following conditions for $\geh$:
\begin{align}
&\;\{i\in I\mid a_{0i}<0\}=\{1\}, \label{eq:ass1}\\
&\;a_{01}=a_{10}=-1.\label{eq:ass2}
\end{align}
Namely, in the Dynkin diagram of $\geh$, $0$ is connected only with $1$ by a single bond.
This implies $\alpha_0=2\La_0-\La_1$. We note that $D_4^{(3)}$ we treat in this paper satisfies
these conditions. We also remark that the labeling of $I$ does not always agree with that of
\cite{Kac}. Let $\gehl{01}$ (resp. $\gehl{\neq0},\gehl{\neq0,1}$) denote the Levi subalgebra of 
$\geh$ corresponding to the index set $\{0,1\}$ (resp. $I\setminus\{0\},I\setminus\{0,1\}$). 
For an integral weight $\lambda$ such that $\langle h_i,\lambda\rangle\ge0$
for $i\not=0$, let us denote
by $B_{\not=0}(\la)$ the $U_q(\geh_{\not=0})$-crystal 
with highest weight $\lambda$.

\begin{theorem} \label{th:uniqueness}
Let $\geh$ be an affine Lie algebra satisfying the conditions \eqref{eq:ass1},
\eqref{eq:ass2}.
Let $B,B'$ be $U'_q(\geh)$-crystals which decompose into 
$\bigoplus_{0\le k\le l}B_{\not=0}(-k\alpha_0)$
as $U_q(\gehl{\neq0})$-crystals. 
Then they are isomorphic to each other as $U'_q(\geh)$-crystals.
\end{theorem}
The theorem says that under the assumptions \eqref{eq:ass1},
\eqref{eq:ass2}, there is a unique 
way to draw $0$-arrows in the crystal graph of the $U_q(\gehl{\neq0})$-crystal 
$\bigoplus_{0\le k\le l}B_{\not=0}(-k\alpha_0)$.

For an element $x$ of the Weyl group $W$ of $\g$,
let $x=s_{i_1}s_{i_2}\cd s_{i_l}$ be a reduced expression of $x$ 
by simple reflections. 
We define $\et{x}^{\max}=
\et{i_1}^{\max}\et{i_2}^{\max}\cd\et{i_l}^{\max}$ and
$S_x=S_{i_1}S_{i_2}\cd S_{i_l}$. Here 
$\et{i}^{\max}b=\et{i}^{\veps_i(b)}b$ and 
$S_i$ is the Weyl group action on crystals. 
Note that
$\et{x}^{\max}$ or $S_x$ do not depend on the choice of a reduced expression. 
For these matters
along with basic notations on crystals, see \cite{K2}.


The rest of this section is devoted 
to the proof of Theorem~\ref{th:uniqueness}. 
Let $B$ and $B'$ be as in Theorem~\ref{th:uniqueness}.
Let $\psi\cl B\to B'$ be a unique $U_q(\geh_{\not=0})$-crystal
isomorphism.
It is enough to show that $\psi$ commutes with $\et{0}$ and $\ft{0}$.

Set $\la=-\alpha_0$. 
For a weight $\mu$ of the form $w(k\lambda)$ ($w\in W$, $0\le k\le l$),
we denote by $u_{\mu}$ a unique element of $B_{\not=0}(k\lambda)$
with weight $\mu$.
We denote by the same letter the corresponding element of $B$.

For $b\in B$ or $b\in B'$, let us denote by $B_{01}(b)$
the connected $U_q(\geh_{01})$-subcrystal containing $b$.
Similarly, we denote by $B_{0}(b)$
the connected $U_q(\geh_{0})$-subcrystal containing $b$.

We prepare several lemmas.
The next lemma is the same as Sublemma 6.2 of \cite{KS}.
Let us denote by $w$ 
the longest element of the Weyl group of $\gehl{\neq0,1}$. 
\begin{lemma} \label{lem:2}
$w\alpha_1=\delta-\alpha_0-\alpha_1=s_1(\delta-\alpha_0)$.
Moreover, the length
$\ell(s_1ws_1ws_1)$ of
$s_1ws_1ws_1$ is equal to $2\ell(w)+3$.
\end{lemma}
Hence, $w\alpha_1=-s_1\alpha_0$ as elements of $P_{cl}$, and
$s_0=s_1ws_1ws_1$ as automorphisms of $P_{cl}$.
Moreover $\et{s_1ws_1ws_1}^{\max}=
\et{1}^{\max}\et{w}^{\max}\et{1}^{\max}\et{w}^{\max}\et{1}^{\max}$.

\begin{lemma} \label{lem:1}
$\et{0}^ku_{l\la}=u_{(l-k)\la}$ for $0\le k\le 2l$.
\end{lemma}
\begin{proof}
One knows that $S_1u_{l\la}$ is the lowest weight vector of 
the $U_q(\gehl{01})$-crystal
$B^{A_2}(l(\La_0+\La_1))$. With the notations
in section \ref{subsec:A2 review}, $\et{0}^ku_{l\la}$ is identified
with $t(0,l,2l-k)$. If $0\le k\le l$, then $\veps_1(\et{0}^ku_{l\la})=0$ from \eqref{A2 eps phi}.
We also have $\veps_i(\et{0}^ku_{l\la})=\veps_i(u_{l\la})=0$ for $i\in I\setminus\{0,1\}$. Hence
$\et{0}^ku_{l\la}$ is a $\g_{\neq0}$-highest vector of weight
$(l-k)\la$ and coincides with $u_{(l-k)\la}$,

Similarly, for $0\le k\le l$,
one can show $\ft{0}^ku_{-l\la}=u_{-(l-k)\la}$, which completes 
the proof.
\end{proof}
By this lemma, $B_{01}(S_1u_{l\lambda})$
contains all the $\g_{\neq0}$-highest weight vectors.

\begin{lemma}\label{lem:1.5}
The restriction of $\psi$ gives a $U_q(\g_{01})$-crystal isomorphism
$$B_{01}(S_1u_{l\alpha_0})\isoto B_{01}(\psi(S_1u_{l\alpha_0})).$$
\end{lemma}
\begin{proof}
Any element of
$B_{01}(S_1u_{l\alpha_0})$ can be written as 
$\ft{1}^d\ft{0}^c\ft{1}^aS_1u_{l\alpha_0}$
with $a\le l$, $a\le c\le a+l$ and $d\le l-2a+c$.
Set $b=\ft{0}^c\ft{1}^aS_1u_{l\alpha_0}$.
We can see easily by section \ref{subsec:A2 review}
\be[(i)]
%
\item 
If $c\ge l$, we have
$b=S_w\ft{1}^{l-a}u_{(a-c)\alpha_0}$.

Indeed, we have $S_wb=\ft{0}^cS_w\ft{1}^aS_1u_{l\alpha_0}$ and
$S_w\ft{1}^aS_1u_{l\alpha_0}$ has weight
$a(\alpha_0+\alpha_1)-l\alpha_1$.
Since the multiplicity of $B$ at this weight is one,
we have $S_w\ft{1}^aS_1u_{l\alpha_0}=\et{0}^{l+a}\et{1}^a
u_{-l(\alpha_0+\alpha_1)}$.
Hence
$S_wb=\et{0}^{l+a-c}\et{1}^au_{-l(\alpha_0+\alpha_1)}=
\et{0}^{l+a-c}\ft{1}^{l-a}u_{-l\alpha_0}=
\ft{1}^{l-a}\et{0}^{l+a-c}u_{-l\alpha_0}=
\ft{1}^{l-a}u_{(a-c)\alpha_0}$.
\label{item:2}
\item
If $a\le c\le l$, we have
$b=S_1S_w\et{1}^{c-a}u_{(l-a)\alpha_0}$.

In this case $b$ is $1$-highest, and
$S_1b=\et{0}^{2l-c}\et{1}^{l+a-c}S_1u_{-l\alpha_0}$.
Hence by applying \eqref{item:2} by reversing the arrows,
we have
$S_1b=S_w\et{1}^{c-a}u_{(l-a)\alpha_0}$.
\ee
\end{proof}

\begin{lemma} \label{lem:4}
\be[{\rm(i)}]
\item
For any $b\in B$,
$\et{s_1ws_1ws_1}^{\max}b$ is $\gehl{\neq0}$-highest.
\item
Assume that $b\in B$ is $\g_{\neq0,1}$-highest and $1$-lowest,
and $S_wb$ is $1$-highest.
Then $\et{s_1w}^{\max}S_1b$ is $\g_{\neq0}$-highest.
\ee
\end{lemma}
\begin{proof}
(i)\quad The claim follows from the fact that $s_1ws_1ws_1\la=-\la$, 
which can be checked by Lemma~\ref{lem:2}.

\noindent
(ii)\quad
$\et{s_1w}^{\max}S_1b=\et{s_1ws_1ws_1}^{\max}S_1S_wb$ is
$\gehl{\neq0}$-highest by (i).
\end{proof}

\begin{lemma} \label{lem:5}
For $k$ such that $0\le k\le l$ and
an element $b$ of $B_{\not=0}(k\la)$, suppose that
$b$ is $1$-highest,
$S_1b$ is $\gehl{\neq0,1}$-highest,
and $S_wS_1b$ is $1$-highest. Then we have $\langle h_0,\wt b\rangle\le-k$.
\end{lemma}
\begin{proof}
By applying the previous lemma for $S_{s_1ws_1}b$, 
one knows that $\et{s_1w}^{\max}b=\et{s_1ws_1ws_1}^{\max} S_{s_1ws_1}b$
is $\gehl{\neq0}$-highest. Hence we have 
\[
\wt(\et{w}^{\max}b)\in s_1(k\la)+\Zn\alpha_1\quad\text{and}\quad
\wt b\in s_1(k\la)+\Zn\alpha_1+\sum_{i\neq0,1}\Z_{\le0}\alpha_i.
\]
Hence we have 
\[
\langle h_0,\wt b\rangle\le\langle h_0,s_1(k\la)\rangle
=\langle h_0+h_1,k\la\rangle=-k.
\]
\end{proof}


Let $A(r)_i$ ($i=1,2,3$) be the following statements:
\begin{align*}
A(r)_1:&\;\text{for $b\in B$ such that
$\Vert\wt b\Vert^2,\,\Vert\wt b+\alpha_0\Vert^2\ge r$,}\\
&\hs{5ex} \text{we have
$\varepsilon_0(b)=\varepsilon_0(\psi(b))$ and $\psi(\et{0}b)=\et{0}\psi(b)$,}
\\
A(r)_2:&\;\text{for $b\in B$ such that
$\Vert\wt b\Vert^2,\,\Vert\wt b-\alpha_0\Vert^2\ge r$,}\\
&\hs{5ex} \text{we have
$\varphi_0(b)=\varphi_0(\psi(b))$ and $\psi(\ft{0}b)=\ft{0}\psi(b)$,}
\\
A(r)_3:&\;
\text{for $b\in B$ such that
$\Vert\wt b\Vert^2\ge r$, we have
$\psi(S_0b)=S_0\psi(b)$.} 
\end{align*}
We prove these statements for any $r\ge0$ 
by the descending induction on $r$. Assume $A(r')$
for $r'>r$.

Then from $A(r')$ we have 

\begin{lemma} \label{lem:6}
For $b\in B$ such that $\Vert\wt b\Vert^2>r$,
there exists a $U_q(\geh_{01})$-crystal isomorphism
$\xi\cl B_{01}(b)\isoto B_{01}(\psi(b))$.
Moreover, we have
$\xi(b')=\psi(b')$ 
for any $b'\in B_{01}(b)$ such that $\Vert\wt b'\Vert^2>r$.
\end{lemma}

The following lemma implies $A(r)_1$ and $A(r)_2$.
\begin{lemma} \label{lem:7}
Assume that $b\in B$ satisfies $\Vert\wt b\Vert^2\ge r$. Then 
there exists a $U_q(\g_0)$-crystal isomorphism $\eta\cl B_0(b)\isoto
B_0(\psi(b))$
such that $\eta(b)=\psi(b)$ and $\eta(b')=\psi(b')$
for any $b'\in B_0(b)$ such that $\Vert\wt b'\Vert^2>r$.
In particular, $\veps_0(b)=\veps_0(\psi(b))$.
\end{lemma}
\begin{proof}

Assume first that
$b$ is not $1$-extremal.
Then $\Vert\wt \et{1}^{\max}b\Vert>\Vert\wt b\Vert$, and hence
by the preceding lemma,
there exists a $U_q(\g_{01})$-isomorphism
$\xi\cl B_{01}(b)\isoto B_{01}(\psi(b))$ such that
$\xi(b')=\psi(b')$ for $b'\in B_{01}(b)$ 
such that $\Vert\wt b'\Vert^2>r$.
Since $\xi(\et{1}^{\max}b)=\psi(\et{1}^{\max}b)$, we have $\xi(b)=\psi(b)$.

Similarly, if $S_wb$ is not $1$-extremal, the assertion holds.

Now we may assume further that $b$ is $\gehl{\neq0,1}$-highest. 
Moreover we nay assume that $b$ and $S_wb$ are $1$-extremal.
If $b$ is $1$-highest, then $b$ is 
$\gehl{\neq0}$-highest and 
the assertion follows from Lemma \ref{lem:1.5}.
Hence one can assume that $b$ is $1$-lowest.
If $S_wb$ is $1$-lowest, then $S_wb$ is 
$\gehl{\neq0}$-lowest and the assertion holds.
Hence one can assume that $S_wb$ is 1-highest.

By Lemma \ref{lem:4}, $\et{s_1w}^{\max} S_1b$ is $\g_{\neq0}$-highest.
Hence, $\et{s_1w}^{\max} S_1b\in B_{01}(S_1u_{l\lambda})$,
and $\et{w}^{\max} S_1b\in B_{01}(S_1u_{l\lambda})$.
Then Lemma \ref{lem:1.5} implies that 
$\zeta=\psi\vert_{B_{01}(S_1u_{l\lambda})}$ gives a $U_q(\g_{01})$-crystal
isomorphism
$$\zeta\cl B_{01}(\et{w}^{\max}S_1b)
\isoto B_{01}(\psi(\et{w}^{\max}S_1b)),$$
which implies that
\eq
&&\varepsilon_0(\et{w}^{\max}S_1b)=
\varepsilon_0(\psi(\et{w}^{\max}S_1b))
\quad\text{and}\quad
\et{0}^{\max}\psi(\et{w}^{\max}S_1b)=
\psi(\et{0}^{\max}\et{w}^{\max}S_1b).
\label{eq:et0psi}
\eneq
%
%
Hence, $\psi$ induces a $U_q(\g_{0})$-crystal
isomorphism
$B_0(\et{w}^{\max}S_1b)\isoto B_0(\psi(\et{w}^{\max}S_1b))$.

Since $\et{w}^{\max}$ commutes with $\et{0}$, $\ft{0}$ and $\psi$,
$\psi$ gives a $U_q(\g_{0})$-crystal
isomorphism
$B_0(S_1b)\isoto B_0(\psi(S_1b))$
which sends $S_1b$ to $\psi(S_1b)$
and $\et{0}^{\max}S_1b$ to $\psi(\et{0}^{\max}S_1b)$.
In particular, 
$\varepsilon_0(S_1b)=\varepsilon_0(\psi(S_1b))$ and
$\psi(\et{0}^{\max}S_1b)=\et{0}^{\max}\psi(S_1b)$.

Since
$S_1b$ is $1$-highest,
$\et{1}^{\max}\et{0}^{\max}S_1b$ is the highest weight element of 
$B_{01}(b)$.
Similarly, $\psi(\et{1}^{\max}\et{0}^{\max}S_1b)=
\et{1}^{\max}\et{0}^{\max}\psi(S_1b)$
is $(0,1)$-highest.
Hence we have a $U_q(\g_{01})$-crystal isomorphism
$\eta\cl B_{01}(b)\isoto B_{01}(\psi(b))$ which sends
$\et{1}^{\max}\et{0}^{\max}S_1b$ to $\psi(\et{1}^{\max}\et{0}^{\max}S_1b)$,
and therefore $\et{0}^{\max}S_1b$ to $\psi(\et{0}^{\max}S_1b)$
and $S_1b$ to $\psi(S_1b)$.
Hence $\eta$ sends $b$ to $\psi(b)$.
If $b'\in B_0(b)\subset B_{01}(b)$ satisfies $\Vert\wt b'\Vert^2>r$,
then we have $\Vert\wt\et{1}^{\max}\et{0}^{\max}S_1b\Vert^2>r$ 
and Lemma~\ref{lem:6}
implies that $\eta(b')=\psi(b')$.
\end{proof}

\begin{lemma} \label{lem:9}
Assume that $b\in B$ satisfies $\Vert\wt b\Vert^2>r$. Then, 
there exists a $U_q(\g_0)$-crystal isomorphism
$\eta\cl B_0(b)\isoto B_0(\psi(b))$
such that $\eta(b')=\psi(b')$
for any $b'\in B_0(b)$ with $\Vert\wt b'\Vert^2\ge r$.
\end{lemma}
\begin{proof}
By Lemma~\ref{lem:7}, there exists a $U_q(\g_0)$-crystal isomorphism
$B_0(b')\isoto B_0(\psi(b'))$
which sends $b'$ to $\psi(b')$ and $b$ to $\psi(b)$.
\end{proof}

\begin{lemma} \label{lem:8}
If $\Vert\wt b\Vert^2\ge r$ and $b$ is not $0$-extremal,
then $\psi(S_0b)=S_0\psi(b)$.
\end{lemma}
\begin{proof}
By Lemma~\ref{lem:7}, there exist
$U_q(\g_0)$-crystal isomorphisms
$\eta\cl B_0(b)\isoto B_0(\psi(b))$
and
$\eta'\cl B_0(S_0b)\isoto B_0(\psi(S_0b))$
such that $\eta(b)=\psi(b)$ and $\eta'(S_0b)=\psi(S_0b)$.
By the assumption,
$\Vert\wt \et{0}^{\max}b\Vert^2>r$.
Hence $\eta(\et{0}^{\max}b)=\psi(\et{0}^{\max}b)=
\psi(\et{0}^{\max}S_0b)=\eta'(\et{0}^{\max}S_0b)$, which implies that
$\eta=\eta'$.
Hence $\psi(S_0b)=\eta'(S_0b)=\eta(S_0b)=S_0\eta(b)=S_0\psi(b)$.
\end{proof}

\Lemma
Assume that $\Vert\wt b\Vert^2\ge r$
and the $(0,1)$-highest weight $\mu$ of $B_{01}(b)$
satisfies $\Vert\mu\Vert^2>r$.
Then there exists a $U_q(\g_{01})$-crystal isomorphism
$\eta\cl B_{01}(b)\isoto B_{01}(\psi(b))$
such that 
$\eta(b')=\psi(b')$ for any $b'\in B_{01}(b)$
such that $\Vert\wt b'\Vert^2\ge r$.
\enlemma
\begin{proof}
Let $b_0\seteq\et{0}^{\max}\et{1}^{\max}\et{0}^{\max}b$
be the $(0,1)$-highest weight vector of
$B_{01}(b)$. Since $\Vert\wt b_0\Vert^2>r$,
Lemma~\ref{lem:6} implies that
there exists a $U_q(\g_{01})$-crystal isomorphism
$\eta\cl B_{01}(b_0)\isoto B_{01}(\psi(b_0))$
such that $\eta(b')=\psi(b')$ for any $b'\in B_{01}(b_0)$
such that $\Vert\wt b'\Vert^2>r$.
In particular $\eta(b_0)=\psi(b_0)$.
It is enough to show that
$\eta(b)=\psi(b)$.

Assume first that $b$ is not $0$-extremal.
Then $\et{0}^{\max}b$ has square length greater than $r$,
and $\eta(\et{0}^{\max}b)=\psi(\et{0}^{\max}b)$.
Then Lemma~\ref{lem:9} implies that $\eta(b)=\psi(b)$.

Hence we may assume that $b$ is $0$-extremal.
Since the case where $b$ is $0$-lowest is similarly proved 
by reversing the arrows, we assume that $b$ is $0$-highest.
Then $b_0=\et{0}^{\max}\et{1}^{\max}b$.
By Lemma~\ref{lem:9}, $\eta(\et{1}^{\max}b)=\psi(\et{1}^{\max}b)$.
Hence, we have  $\eta(b)=\psi(b)$.
\end{proof}

Now we are ready to complete the proof of $A(r)_3$.
\Lemma
If $\Vert\wt b\Vert^2\ge r$,
then $\psi(S_0b)=S_0\psi(b)$.
\enlemma
\begin{proof}
We shall argue by the descending induction 
on the length of the $0$-string containing $b$.
If $b$ is not $(0,1)$-extremal, then
the preceding lemma implies the desired result.
Hence we may assume that $b$ is $(0,1)$-extremal.

Since the case when $b$ is $0$-lowest is similar,
we assume that $b$ is $0$-highest.
One can assume that $b$ is $\gehl{\neq0,1}$-highest. 
If $b$ is $1$-highest, 
it is $\gehl{\neq0}$-highest and 
the assertion follows from Lemma~\ref{lem:1.5}.
Hence one can assume that $b$ is $1$-lowest.
Similarly, one can assume that $S_wb$ is 1-highest. 

We divide the proof into two cases:
(1) $S_1b$ is $0$-highest, (2) $S_1b$ is not $0$-highest but is $0$-lowest.

First we consider the case (1).
Since $S_1b$ is $(0,1)$-highest in this case, 
we have $\langle h_0,\wt S_1b\rangle\ge0$.
On the other hand, Lemma \ref{lem:5} 
implies that $\langle h_0,\wt S_1b\rangle\le-k$ with $k\ge0$ 
defined by $S_1b\in B_{\neq0}(k\la)$. 
Hence we obtain $k=0$, and $\wt b=0$, which implies
that $S_0b=b$, $S_0\psi(b)=\psi(b)$.

Next, we consider the case (2). In this case, $S_0b$ is 
$(0,1)$-lowest and $S_0S_1b$ is 
$(0,1)$-highest. 
Since $\et{1}^{\max}\et{w}^{\max}S_1b$ is $g_{\neq0}$-highest 
by Lemma~\ref{lem:4},
Lemma~\ref{lem:1.5} implies that
$\psi(S_0\et{w}^{\max}S_1b)=
S_0\psi(\et{w}^{\max}S_1b)$.
Write $\et{w}^{\max}S_1b=\et{i_1}\cdots\et{i_m}S_1b$
with $i_1,\ldots,i_m\in I\setminus\{0,1\}$.
Since $\et{i_\nu}$ commutes with $S_0$ and $\psi$,
we have
$\et{i_1}\cdots\et{i_m}\psi(S_0S_1b)=
\psi(S_0\et{w}^{\max}S_1b)=S_0\psi(\et{w}^{\max}S_1b)
=\et{i_1}\cdots\et{i_m}S_0\psi(S_1b)$.
Hence we have
\eq\label{eq:s0s1}
&&\psi(S_0S_1b)=
S_0\psi(S_1b)=S_0S_1\psi(b).
\eneq
Let $\mu$ be the weight of
$S_0S_1b$, and set $\mu_i=\langle h_i,\mu\rangle$. We have $\mu_0>0$,
for, otherwise, $S_1b$ is $0$-highest.

Note that $\vphi_0(S_0S_1S_0b)=\mu_0+\mu_1>\vphi_0(b)=\mu_1$.
By the descending induction on 
the length of the 0-string containing $b$, we see
that $\psi(S_0S_1S_0b)=S_0\psi(S_1S_0b)$.
Hence
$S_0S_1\psi(S_0b)=\psi(S_0S_1S_0b)=
\psi(S_1S_0S_1b)=S_1\psi(S_0S_1b)=S_1S_0S_1\psi(b)
=S_0S_1S_0\psi(b)$.
Here, we used \eqref{eq:s0s1}.
Hence we obtain $\psi(S_0b)=S_0\psi(b)$.
\end{proof}

Thus, the descending induction on $r$ proceeds,
and the proof of Theorem~\ref{th:uniqueness} is complete.

\section{Perfectness of the crystal $B_l$}

\subsection{Connectedness of $B_l\ot B_l$ (Proof of (P1))}

We show that any element $b\ot b'$ of $B_l\otimes B_l$ can be connected with 
$\phi\otimes\phi$. Here $\phi$ stands for $(0,0,0,0,0,0)\in B_l$. Like this we use in this section
the tableau representation for an element of $B_l$. By applying $\ft{1}$ and $\ft{2}$ sufficiently 
many times, one can assume that $\ft{i}(b\ot b')=0$ ($i=1,2$). This implies $\ft{i}b'=0$ ($i=1,2$),
namely, 
\[
b'={\bar 1}^m\text{ for some }0\le m\le l.
\]
Note that the 0-string containing ${\bar 1}^m$ is given by 
\[
{\bar 1}^l\rightarrow\cd\rightarrow{\bar 1}^m\rightarrow\cd\rightarrow{\bar 1}
\rightarrow\phi\rightarrow1\rightarrow\cd\rightarrow1^l.
\]
Set $\gamma(b)=m+(\varphi_0(b)-l+m)_+$, then from the tensor product rule of crystals we have
\[
{\tilde f}_0^{\gamma(b)}(b\otimes {\bar 1}^m)=
{\tilde f}_0^{\gamma(b)-m}(b)\otimes {\tilde f}_0^m({\bar 1}^m)=
{\tilde b}\otimes\phi, 
\quad \text{where ${\tilde b}={\tilde f}_0^{\gamma(b)-m}b$}. 
\]
Since $\ft{i}\phi=0$ for $i=1,2$ there exists a sequence 
$\{i_1,\ldots,i_k\}\subset \{1,2\}$ and a non-negative integer $m'$ such that 
${\tilde f}_{i_k}\ldots{\tilde f}_{i_1}({\tilde b}\otimes\phi)={\bar 1}^{m'}\otimes\phi$.
Thus we have 
\[
{\tilde f}_0^{m'}{\tilde f}_{i_k}\ldots{\tilde f}_{i_1}{\tilde f}_0^{\gamma(b)}
(b\otimes{\bar 1}^m)=\phi\otimes\phi. 
\]

\subsection{Minimal elements in $B_l$ (Proof of (P4) and (P5))}

First we are to show $\langle c,\vphi(b)\rangle\ge l$ for $b\in B_l$. 
{}From Proposition \ref{prop:eps0-phi0} and formulas of $\veps_i,\vphi_i$ 
($i=1,2$) in section \ref{subsec:def} we have
\begin{align*}
\langle c, \varphi(b)\rangle
=&\varphi_0(b)+2\varphi_1(b)+3\varphi_2(b) \\
=&l+\max A+2(z_3+(z_2)_+)_++(3z_4)_+ -(z_1+z_2+2z_3+3z_4),
\end{align*}
where $z_j$ ($1\le j\le 4$) are given in \eqref{z1-4} and $A$ is given in \eqref{A}.

\begin{lemma}\label{lem:PerfectCrystal}
For $(z_1,z_2,z_3,z_4)\in \mathbb{Z}^4$ set 
\[
\psi(z_1,z_2,z_3,z_4)=\max A+
2(z_3+(z_2)_+)_+ +(3z_4)_+ -(z_1+z_2+2z_3+3z_4).
\]
Then we have $\psi(z_1,z_2,z_3,z_4)\geq 0$ and $\psi(z_1,z_2,z_3,z_4)=0$ if and only if 
$(z_1,z_2,z_3,z_4)=(0,0,0,0)$. 
\end{lemma}

\begin{proof}
Note that for $z\in\Z$, $2(z)_+\ge z$ and $2(z)_+=z$ implies $z=0$. 
We prove by dividing the cases of the values that attain the maximum in $\max A$.

Suppose $\max A=0$. Using the above inequality, we have 
\[
\psi\ge (z_2)_++(3z_4)_+-(z_1+z_2+z_3+3z_4). 
\]
Since $z_1+z_2+z_3+3z_4\le0$, we have $\psi\ge0$.
$\psi=0$ holds if and only if 
\[
z_1+z_2+3z_4=0,\;z_3=0,\;z_2,z_4\le0.
\]
Since we have $z_1\le0$ in this case, one can conclude that $\psi=0$ implies 
$(z_1,z_2,z_3,z_4)=(0,0,0,0)$.

The other cases are similar. In particular, if $\max A>0$, $\psi>0$.
\end{proof}

Thanks to the lemma, we have 
\[
\langle c,\varphi(b)\rangle-l = \psi(z_1,z_2,z_3,z_4)\geq 0.
\]
Since $\langle c,\vphi(b)-\veps(b)\rangle=0$, we also obtain $\langle c,\veps(b)\rangle\ge l$,
which proves (P4).

Suppose $\langle c,\veps(b)\rangle=l$. It implies $\psi=0$. Hence from the lemma one can conclude 
that such element $b=(x_1,x_2,x_3,\bar{x}_3,\bar{x}_2,\bar{x}_1)$ should satisfy
$x_1={\bar x}_1,\;x_2=x_3={\bar x}_3={\bar x}_2$. Therefore we have 
\[
(B_l)_{\min}=\{(\alpha,\beta,\beta,\beta,\beta,\alpha)
\vert \alpha, \beta \in \Zn,2\alpha+3\beta\le l\}. 
\]
For $b=(\alpha,\beta,\beta,\beta,\beta,\alpha)\in B_l$ one calculates 
\[
\veps(b)=\vphi(b)=(l-2\alpha-3\beta)\La_0+\alpha\La_1+\beta\La_2. 
\]
Thus we have also shown (P5).

\subsection{Coherent family of perfect crystals}

We review the notion of a coherent family of perfect crystals introduced in \cite{KKM}. 
Let $\{B_l\}_{l\ge1}$ be a family of perfect crystals $B_l$ of level $l$ and $(B_l)_{\min}$ be
the subset of minimal elements of $B_l$. Set $J=\{(l,b)\mid l\in\Z_{>0},b\in(B_l)_{\min}\}$.
Let $\sigma$ denote the isomorphism of $(P^+_{cl})_l$ defined by $\sigma=\veps\circ\vphi^{-1}$.

\begin{define} \label{def:coherent}
A crystal $B_\infty$ with an element $b_\infty$ is called a limit of $\{B_l\}_{l\ge1}$ if it 
satisfies the following conditions:
\begin{itemize}
\item[$\bullet$] $\wt b_\infty=0,\veps(b_\infty)=\vphi(b_\infty)=0$,

\item[$\bullet$] for any $(l,b)\in J$, there exists an embedding of crystals
\[
f_{(l,b)}:\;T_{\veps(b)}\ot B_l\ot T_{-\vphi(b)}\longrightarrow B_\infty
\]
sending $t_{\veps(b)}\ot b\ot t_{-\vphi(b)}$ to $b_\infty$,

\item[$\bullet$] $B_\infty=\bigcup_{(l,b)\in J}\mbox{Im}\,f_{(l,b)}$.
\end{itemize}
If a limit exists for the family $\{B_l\}$, we say that $\{B_l\}$ is a coherent family of perfect
crystals.
\end{define}
For $\la\in P_{cl}$ $T_{\la}$ denotes a crystal with a unique element $t_{\la}$. See \cite{K2} for 
the details. For our purpose the following facts are sufficient. For any $P_{cl}$-weighted crystal
$B$ and $\la,\mu\in P_{cl}$ consider the crystal
\[
T_{\la}\ot B\ot T_{\mu}=\{t_{\la}\ot b\ot t_{\mu}\mid b\in B\}.
\]
The crystal structure is given by
\begin{align*}
\et{i}(t_{\la}\ot b\ot t_{\mu})&=t_{\la}\ot\et{i}b\ot t_{\mu}, & 
\ft{i}(t_{\la}\ot b\ot t_{\mu})&=t_{\la}\ot\ft{i}b\ot t_{\mu}, \\
\veps_i(t_{\la}\ot b\ot t_{\mu})&=\veps_i(b)-\langle h_i,\la\rangle, &
\vphi_i(t_{\la}\ot b\ot t_{\mu})&=\vphi_i(b)+\langle h_i,\mu\rangle, \\
\wt(t_{\la}\ot b\ot t_{\mu})&=\la+\mu+\wt b.
\end{align*}

Let us now consider the following set
\[
B_\infty=\{b=(\nu_1,\nu_2,\nu_3,\bar{\nu}_3,\bar{\nu}_2,\bar{\nu}_1)\mid
\nu_i,\bar{\nu}_i\in\Z,\nu_3\equiv\bar{\nu}_3\;(\text{mod }2)\},
\]
and set $b_\infty=(0,0,0,0,0,0)$. We introduce the crystal structure on $B_\infty$ as follows.
The actions of $\et{i},\ft{i}$ ($i=0,1,2$) are defined by the same rule as in section \ref{subsec:def}
and \ref{subsec:0-action} with $x_i$ and $\bar{x}_i$ replaced with $\nu_i$ and $\bar{\nu}_i$. The 
only difference lies in the fact that $\et{i}b$ or $\ft{i}$ never become $0$, since we allow a 
coordinate to be negative and there is no restriction for the sum $s(b)=\sum_{i=1}^2(\nu_i+\bar{\nu}_i)
+(\nu_3+\bar{\nu}_3)/2$. For $\veps_i,\vphi_i$ with $i=1,2$ we adopt the formulas in section 
\ref{subsec:def}. For $\veps_0,\vphi_0$ we define
\begin{align*}
\veps_0(b)&=-s(b)+\max A-(2z_1+z_2+z_3+3z_4),\\
\vphi_0(b)&=-s(b)+\max A,
\end{align*}
where 
\[
A=(0,z_1,z_1+z_2,z_1+z_2+3z_4,z_1+z_2+z_3+3z_4,2z_1+z_2+z_3+3z_4)
\]
and $z_1,z_2,z_3,z_4$ are given in \eqref{z1-4} with $x_i,\bar{x}_i$ replaced with $\nu_i,\bar{\nu}_i$.
Note that $\wt b_\infty=0$ and $\veps_i(b_\infty)=\vphi_i(b_\infty)=0$ for $i=0,1,2$.

Let $b_0=(\alpha,\beta,\beta,\beta,\beta,\alpha)$ be an element of $(B_l)_{\min}$. Since 
$\veps(b_0)=\vphi(b_0)$, one can set $\sigma=\mbox{id}$. Let $\la=\veps(b_0)$. For 
$b=(x_1,x_2,x_3,\bar{x}_3,\bar{x}_2,\bar{x}_1)\in B_l$ we define a map
\[
f_{(l,b_0)}:\;T_{\la}\ot B_l\ot B_{-\la}\longrightarrow B_\infty
\]
by
\[
f_{(l,b_0)}(t_{\la}\ot b\ot t_{-\la})=b'=(\nu_1,\nu_2,\nu_3,\bar{\nu}_3,\bar{\nu}_2,\bar{\nu}_1)
\]
where
\begin{align*}
\nu_1&=x_1-\alpha, & \bar{\nu}_1&=\bar{x}_1-\alpha, \\
\nu_j&=x_j-\beta, & \bar{\nu}_j&=\bar{x}_j-\beta\;(j=2,3).
\end{align*}
Then we have 
\begin{align*}
\wt(t_{\la}\ot b\ot t_{-\la})&=\wt b=\wt b', \\
\vphi_0(t_{\la}\ot b\ot t_{-\la})&=\vphi_0(b)+\langle h_0,-\la\rangle \\
&\vphi_0(b')+(l-s(b))+s(b')-(l-2\alpha-3\beta)=\vphi_0(b'),\\
\vphi_1(t_{\la}\ot b\ot t_{-\la})&=\vphi_1(b)+\langle h_1,-\la\rangle=\vphi_1(b')+\alpha-\alpha
=\vphi_1(b'), \\
\vphi_2(t_{\la}\ot b\ot t_{-\la})&=\vphi_2(b)+\langle h_2,-\la\rangle=\vphi_2(b')+\beta-\beta
=\vphi_2(b').
\end{align*}
$\veps_i(t_{\la}\ot b\ot t_{-\la})=\veps_i(b')$ ($i=0,1,2$) also follows from the above calculations.

It is straightforward to check that if $b,\et{i}b\in B_l$ (resp. $b,\ft{i}b\in B_l$), then 
$f_{(l,b_0)}(\et{i}(t_{\la}\ot b\ot t_{-\la}))=\et{i}f_{(l,b_0)}(t_{\la}\ot b\ot t_{-\la})$
(resp. $f_{(l,b_0)}(\ft{i}(t_{\la}\ot b\ot t_{-\la}))=\ft{i}f_{(l,b_0)}(t_{\la}\ot b\ot t_{-\la})$).
Hence $f_{(l,b_0)}$ is a crystal embedding. It is easy to see that 
$f_{(l,b_0)}(t_{\la}\ot b_0\ot t_{-\la})=b_\infty$. We can also check 
$B_\infty=\bigcup_{(l,b)\in J}\mbox{Im}\,f_{(l,b)}$. Therefore we have shown that the family of
perfect crystals $\{B_l\}_{l\ge1}$ forms a coherent family.

\section*{Acknowledgments}
\smallskip\par\noindent
KCM thanks the staff of Osaka University for their hospitality during his visit in 2004
and acknowledges partial support from NSA grant H98230-06-1-0025.
MK and MO are partially supported by
Grant-in-Aid for Scientific Research (B) 18340007 and (C) 18540030,
Japan Society for the Promotion of Science, respectively.
DY is supported by the 21 century COE program
at Graduate School of Mathematical Sciences, the University of Tokyo.

\appendix

\section{Table of $\ft{1}^q\ft{0}^p{\bar b}_{j_0,j_1}^{l,i}$} \label{app:A}

Assume $j_0\le j_1$. We give the table of $x=\ft{1}^q\ft{0}^p{\bar b}_{j_0,j_1}^{l,i}$ on $B_{\ge0}$ 
for $0\le p\le j_0$, $0\le q\le j_1+p$.

\begin{itemize}

\item[I.] $0\leq p\leq i$ case: 

\begin{itemize}

\item[(i)] $0\leq q\leq j_0-i+p$ case: 
\[
x=(p,y_1,3y_0-2y_1+i-p,y_0+i-p,y_0+j_0-q,q).
\]

\item[(ii)] $j_0-i+p\leq q\leq j_1$ case: 
\[
x=(p,y_1,3y_0-2y_1-q+j_0,y_0+2i-2p+q-j_0,y_0+i-p,j_0-i+p).
\]

\item[(iii)] $j_1\leq q\leq j_1+p$ case: 
\[
x=(p-q+j_1,y_1+q-j_1,y_1,y_0+2i-2p+j_1-j_0,y_0+i-p,j_0-i+p).
\]

\end{itemize}

\item[II.] $i\le p\le j_0$ case: 

\begin{itemize}

\item[(i)] $0\leq q\leq j_0-p+i$ case: 
\[
x=(i,y_1,2p-2i+3y_0-2y_1,y_0,y_0+j_0-p+i-q,q).
\]

\item[(ii)] $j_0-p+i\leq q\leq p-i+j_1$ case: 
\[
x=(i,y_1,3y_0-2y_1-i+j_0+p-q,y_0+q-j_0+p-i,y_0,j_0-p+i).
\]

\item[(iii)] $p-i+j_1\leq q\leq j_1+p$ case: 
\[
x=(p+j_1-q,y_1+q-p+i-j_1,y_1,2p-2i+j_1-j_0+y_0,y_0,j_0-p+i).
\]

\end{itemize}
\end{itemize}

\section{Table of $\ft{0}^r\ft{1}^q\ft{0}^p{\bar b}_{i,j_1}^{l,i}$} \label{app:B}

We give the table of $x=\ft{0}^r\ft{1}^q\ft{0}^p{\bar b}_{i,j_1}^{l,i}$ on $B_{\ge0}$ for 
$0\le p\le q\le j_1+p,0\le r\le i+q-2p$.

\begin{itemize}

\item[I.] $p\leq i$, $p\le q$, $r\leq i+q-2p$ case: 

\begin{itemize}

\item[(i)] $q\le p+\frac{1}{2}(j_1-i)$, $0\leq r\leq i+q-2p$ case: 
\[
x=(p+r,y_1,j_1+y_1-q-r,y_0+i-2p+q-r,y_0+i-p,p).
\]

\item[(ii)] $p+\frac{1}{2}(j_1-i)\leq q\le j_1$, $0\leq r\le j_1-q$ case: same as (i).

\item[(iii)] $p+\frac{1}{2}(j_1-i)\leq q\le p+\frac{2}{3}(j_1-i)$, $j_1-q\leq r\leq i+q-2p$ case: 
\[
x=(j_1+p-q,y_1-j_1+q+r,y_1,y_0+j_1+i-2p-2r,y_0+i-p,p).
\]

\item[(iv)] $p+\frac{2}{3}(j_1-i)\leq q\le j_1$, $j_1-q\leq r\le i-p+\frac{2}{3}(j_1-i)$ case: 
same as (iii).

\item[(v)] $p+\frac{2}{3}(j_1-i)\leq q\le j_1$, $i-p+\frac{2}{3}(j_1-i)\leq r\leq i+q-2p$ case: 
\[
x=(j_1+p-q,2y_1-y_0-p+q,2y_0-y_1-2j_1+2p+2r,y_1,y_1+j_1+i-2p-r,p).
\]

\item[(vi)] $j_1\leq q\leq j_1+p$, $0\leq r\leq i-p+\frac{2}{3}(j_1-i)$ case: same as (iii).

\item[(vii)] $j_1\leq q\leq j_1+p$, $i-p+\frac{2}{3}(j_1-i)\le r\leq j_1+i-2p$ case: 
\[
x=(j_1+p-q,2y_1-y_0-p+q,2y_1-y_0-j_1-i+2p+2r,y_1,y_1+j_1+i-2p-r,p).
\]

\item[(viii)] $j_1\leq q\leq j_1+p$, $j_1+i-2p\leq r\leq i+q-2p$ case: 
\[
x=(j_1+p-q,2y_1-y_0-p+q,2y_1-y_0+r,y_1-j_1-i+2p+r,y_1,j_1+i-p-r).
\]

\end{itemize}

\item[II.] $p\leq i$, $p\le q$, $r\ge i+q-2p$ case: 

\begin{itemize}

\item[(i)] $p\le q\leq p+\frac{1}{2}(j_1-i)$ case: 
\[
x=(p+r,y_1,y_1+j_1-i+2p-2q,y_0,y_0+i-p,p).
\]

\item[(ii)] $p+\frac{1}{2}(j_1-i)\le q\leq p+\frac{2}{3}(j_1-i)$ case: 
\[
x=(j_1-i+3p-2q+r,y_1-j_1+i-2p+2q,y_1,y_0+j_1-i+2p-2q,y_0+i-p,p).
\]

\item[(iii)] $p+\frac{2}{3}(j_1-i)\le q\le j_1$ case: 
\[
x=(j_1-i+3p-2q+r,2y_1-y_0-p+q,2y_0-y_1-2j_1+2i-2p+2q,y_1,y_1+j_1-q,p).
\]

\item[(iv)] $j_1\leq q\leq j_1+p$ case: 
\[
x=(j_1-i+3p-2q+r,2y_1-y_0-p+q,2y_1-y_0+i+q-2p,y_1-j_1+q,y_1,j_1+p-q).
\]

\end{itemize}
\end{itemize}

\section{Table of $\ft{0}^r\ft{1}^q\ft{0}^{j_0}{\bar b}_{j_0,j_1}^{l,i}$} \label{app:C}

Assume $j_0\le j_1$.
We give the table of $x=\ft{0}^r\ft{1}^q\ft{0}^{j_0}{\bar b}_{j_0,j_1}^{l,i}$ on $B_{\ge0}$ for 
$0\le q\le j_0+j_1,0\le r$.

\begin{itemize}

\item[I.] $0\le q\le i$ case:
\[
x=(i+r,y_1,y_1+j_0+j_1-2i,y_0,y_0+i-q,q).
\]

\item[II.] $i\le q\le j_0$ case:
\[
x=(i+r,y_1,y_1+j_0+j_1-i-q,y_0-i+q,y_0,i).
\]

\item[III.] $j_0\le q\le j_0+j_1-i$ case:

\begin{itemize}

\item[(i)] $j_0\le q\le j_0+\frac{j_1-i}2,0\le r\le q-j_0$ or 
$j_0+\frac{j_1-i}2\le q\le j_0+j_1-i,0\le r\le j_0+j_1-i-q$ case:
\[
x=(i+r,y_1,y_1+j_0+j_1-i-q-r,y_0-i+q-r,y_0,i).
\]

\item[(ii)] $j_0\le q\le j_0+\frac{j_1-i}2,r\ge q-j_0$ case:
\[
x=(i+r,y_1,y_1+2j_0+j_1-i-2q,y_0+j_0-i,y_0,i).
\]

\item[(iii)] $j_0+\frac{j_1-i}2\le q\le \frac{4j_0+2j_1}3-i,j_0+j_1-i-q\le r\le q-j_0$ or
$\frac{4j_0+2j_1}3-i\le q\le j_0+j_1-i,j_0+j_1-i-q\le r\le \frac{j_0+2j_1}3-i$ case:
\[
x=(j_0+j_1-q,y_1-j_0-j_1+i+q+r,y_1,y_0+j_0+j_1-2i-2r,y_0,i).
\]

\item[(iv)] $j_0+\frac{j_1-i}2\le q\le \frac{4j_0+2j_1}3-i,r\ge q-j_0$ case:
\[
x=(2j_0+j_1-2q+r,y_1-2j_0-j_1+i+2q,y_1,y_0+3j_0+j_1-2i-2q,y_0,i).
\]

\item[(v)] $\frac{4j_0+2j_1}3-i\le q\le j_0+j_1-i,\frac{j_0+2j_1}3-i\le r\le q-j_0$ case:
\[
x=(j_0+j_1-q,2y_1-y_0-j_0+q,2y_0-y_1-2j_1+2i+2r,y_1,y_1+j_1-i-r,i).
\]

\item[(vi)] $\frac{4j_0+2j_1}3-i\le q\le j_0+j_1-i,r\ge q-j_0$ case:
\[
x=(2j_0+j_1-2q+r,2y_1-y_0-j_0+q,
2y_0-y_1-2j_0-2j_1+2i+2q,y_1,y_1+j_0+j_1-i-q,i).
\]

\end{itemize}

\item[IV.] $j_0+j_1-i\le q\le j_0+j_1$ case:

\begin{itemize}

\item[(i)] $0\le r\le \frac{j_0+2j_1}3-i$ case: same as III (iii).

\item[(ii)] $\frac{j_0+2j_1}3-i\le r\le j_1-i$ case: same as III (v).

\item[(iii)] $j_1-i\le r\le q-j_0$ case:
\[
x=(j_0+j_1-q,2y_1-y_0-j_0+q,2y_0-y_1-j_1+i+r,y_1-j_1+i+r,y_1,j_1-r).
\]

\item[(iv)] $r\ge q-j_0$ case:
\[
x=(2j_0+j_1-2q+r,2y_1-y_0-j_0+q,2y_0-y_1-j_0-j_1+i+q,y_1-j_0-j_1+i+q,y_1,j_0+j_1-q).
\]

\end{itemize}
\end{itemize}

\section{Table of $\ft{0}^{j_0+q-2p}\ft{1}^q\ft{0}^p{\bar b}_{j_0,j_1}^{l,i}$} \label{app:D}

Assume $j_0\le j_1$.
We give the table of $x=\ft{0}^{j_0+q-2p}\ft{1}^q\ft{0}^p{\bar b}_{j_0,j_1}^{l,i}$ on $B_{\ge0}$ for 
$0\le p\le j_0,p\le q\le j_1+p$.

\begin{itemize}

\item[I.] $0\le p\le i$ case:

\begin{itemize}

\item[(i)] $i\le q-p+i\le \frac{j_0+j_1}2$ case:
\[
x=(i-p+q,y_1,y_1+j_0+j_1-2i+2p-2q,y_0,y_0+i-p,p).
\]

\item[(ii)] $\frac{j_0+j_1}2\le q-p+i\le \frac{j_0+2j_1}3$ case:
\[
x=(j_0+j_1-i+p-q,y_1-j_0-j_1+2i-2p+2q,y_1,y_0+j_0+j_1-2i+2p-2q,y_0+i-p,p).
\]

\item[(iii)] $q-p+i\ge \frac{j_0+2j_1}3,q\le j_1$ case:
\[
x=(j_0+j_1-i+p-q,2y_1-y_0-j_0+i-p+q,2y_0-y_1-2j_1+2i-2p+2q,y_1,y_1+j_1-q,p).
\]

\item[(iv)] $j_1\le q\le j_1+p$ case:
\[
x=(j_0+j_1-i+p-q,2y_1-y_0-j_0+i-p+q,2y_1-y_0-j_0+2i-2p+q,y_1-j_1+q,y_1,j_1+p-q).
\]

\end{itemize}

\item[II.] $i\le p\le j_0$ case:

\begin{itemize}

\item[(i)] $2p\le 2q\le j_0+j_1+p-i$ case:
\[
x=(i-p+q,y_1,y_1+j_0+j_1-i+p-2q,y_0-i+p,y_0,i).
\]

\item[(ii)] $2q\ge j_0+j_1+p-i,q-p+i\le \frac{j_0+2j_1}3$ case:
\[
x=(j_0+j_1-q,y_1-j_0-j_1+i-p+2q,y_1,y_0+j_0+j_1-2i+2p-2q,y_0,i).
\]

\item[(iii)] $\frac{j_0+2j_1}3\le q-p+i\le j_1$ case:
\[
x=(j_0+j_1-q,2y_1-y_0-j_0+q,2y_0-y_1-2j_1+2i-2p+2q,y_1,y_1+j_1-i+p-q,i).
\]

\item[(iv)] $j_1\le q-p+i\le j_1+i$ case:
\[
x=(j_0+j_1-q,2y_1-y_0-j_0+q,2y_0-y_1-j_1+i-p+q,y_1-j_1+i-p+q,y_1,j_1+p-q).
\]

\end{itemize}
\end{itemize}

\end{document}